\documentclass[a4paper,11pt]{scrartcl}

\usepackage[USenglish]{babel}
\usepackage[T1]{fontenc}
\usepackage[utf8]{inputenc}

\usepackage{graphicx,subfigure,tikz}

\usepackage{amsmath,amssymb,amsfonts,amsthm,array,stmaryrd,booktabs,paralist,todonotes,mathtools}

\usepackage{ifpdf}
\ifpdf
\usepackage[pdftex,colorlinks=true,linkcolor=blue,citecolor=green,urlcolor=blue,bookmarks]{hyperref}
\else
\usepackage[dvipdfm,colorlinks=true,linkcolor=blue,citecolor=green,urlcolor=blue,bookmarks]{hyperref}
\fi 

\newcounter{para}[section]
\numberwithin{equation}{section}

\newcommand{\dd}{\textup{d}}

\newtheorem{thm}{Theorem}[section]
\newtheorem{lem}[thm]{Lemma}
\newtheorem{defi}[thm]{Definition}
\newtheorem{rem}[thm]{Remark}
\newtheorem{cor}[thm]{Corollary}
\newenvironment{prf}{\textbf{Proof:}} {\hspace*{\fill} $\square$ \newline}

\begin{document}

    \title{A general existence result for isothermal two-phase flows with phase transition}
    \author{Maren Hantke\footnotemark[1]\,\, and Ferdinand Thein\footnotemark[2] \footnotemark[3]}
    \date{\today}
    \maketitle
    \begin{abstract}
        Liquid-vapor flows with phase transitions have a wide range of applications.
        Isothermal two-phase flows described by a single set of isothermal Euler equations,
        where the mass transfer is modeled by a kinetic relation, have been investigated analytically in (Quarterly of applied Mathematics, vol.\ LXXI 3 (2013), pp.\ 509-540.).
        This work was restricted to liquid water and its vapor modeled by linear equations of state.
        The focus of the present work lies on the generalization of the primary results to arbitrary substances, arbitrary equations of state and thus a more general kinetic relation.
        We prove existence and uniqueness results for Riemann problems. In particular, nucleation and evaporation are discussed.
    \end{abstract}
    \renewcommand{\thefootnote}{\fnsymbol{footnote}}
    \footnotetext[1]{Institute for Analysis and Numerics, Otto-von-Guericke University Magdeburg,\\ PSF 4120, D-39016 Magdeburg, Germany.
    \href{mailto:maren.hantke@ovgu.de}{\textit{maren.hantke@ovgu.de}}}
    \footnotetext[2]{Institute for Analysis and Numerics, Otto-von-Guericke University Magdeburg,\\ PSF 4120, D-39016 Magdeburg, Germany.
    \href{mailto:ferdinand.thein@ovgu.de}{\textit{ferdinand.thein@ovgu.de}}}
    \footnotetext[3]{The second author is supported by the DFG grant HA-6471/2-1. The authors thankfully acknowledge the support.}
    \renewcommand{\thefootnote}{\arabic{footnote}}
    \newpage

    \tableofcontents
    \newpage

    \section{Introduction}\label{chap_intro}
\setcounter{para}{1}
Compressible liquid-vapor flows have a wide range of applications. Two-phase flow models are 
used to describe such processes, e.g.\ the formation of clouds, cavitation near moving objects in liquids such as ship propellers or certain phenomena in biology.
Main difficulties in the modeling result from the phase interactions, especially from mass and energy transfer due to condensation or evaporation processes.
Several two-phase flow models are available in the literature. They are mainly distinguished in sharp and diffusive interface models.
For a detailed discussion of these models we refer to Zein \cite{Zein2010a} and concerning sharp interface models we exemplary refer to Bedeaux et al. \cite{Bedeaux2004}.
In our work we study compressible two-phase flows with phase transitions across a sharp interface. Phase transitions are modeled using a kinetic relation.
This concept was introduced by Abeyaratne and Knowles \cite{Abeyaratne1991} for solid-solid phase transitions.
This kinetic relation controls the mass transfer across the interface between the two adjacent phases.
For a more general context of kinetic relations see LeFloch \cite{LeFloch2002}.
A detailed and very interesting survey on the Riemann problem for a large class of thermodynamic consistent constitutive models in the setting of Euler equations models
can be found in Menikoff and Plohr \cite{Menikoff1989}. Here the considerations are restricted to a simple kinetic relation that results from the assumption of local equilibrium at the interface.\\
In a recent work by Hantke et al.\ \cite{Hantke2013} Riemann problems relying on the isothermal Euler equations with a non-monotone pressure-density function are considered.
This function is composed of three parts: the equations of state for the two single phases and an arbitrary relation for the intermediate state.
The two phases are distinguished using the Maxwell construction, also known as the Equal-Area-Rule. The mass transfer is modeled via a kinetic relation, derived in \cite{Dreyer2012},
based on classical Hertz-Knudsen theory, see \cite{Bond2004}.
The authors discussed Riemann problems for various different cases of initial data and showed existence and uniqueness. Furthermore Hantke et al.\ also covered the cases of cavitation and nucleation.
The constructed Riemann solutions are selfsimilar. They consist of constant states, separated by classical rarefaction and shock waves or phase boundaries.
Nevertheless, the basic assumptions are very restrictive. Existence and uniqueness results are proven for liquid water and its vapor, modeled by linear equations of state.\\
Also M\"uller and Voss \cite{Mueller2006}, \cite{Voss2005} considered the isothermal Euler system. In contrast to the above mentioned work they modeled the fluid
using the van der Waals equation of state. Instead of a kinetic relation the Liu entropy condition is used to achieve uniqueness.
As a consequence M\"uller and Voss need non-classical composite waves to construct solutions.
Further literature in this context is given by Merkle \cite{Merkle2006}, Merkle and Rohde \cite{Merkle2007}.
The focus of our present work is on the distinguished generalization of the results of Hantke et al.\ \cite{Hantke2013} resp.\  Menikoff and Plohr \cite{Menikoff1989}.
We consider two-phase flows for any regular fluid. Both phases can be modeled by any thermodynamic relevant equation of state.
Further we construct exact Riemann solutions and prove existence and uniqueness results that advance achievements in the actual literature.\\
The paper is organized as follows. In Section \ref{chap_euler_eq} we present the balance equations in the bulk phases and the corresponding jump conditions across discontinuities.
Further we give the thermodynamic framework needed throughout this work and discuss the Riemann problem in the isothermal case including the entropy inequality.
In Section \ref{chap_pb_sol} we prove existence and uniqueness of a solution at the interface under certain appropriate assumptions.
The following Section \ref{chap_2p_rp} contains a monotonicity argument needed to solve the two-phase Riemann problem, which is done subsequently.
In Section \ref{chap_3p_rp} we present solutions to initial one-phase Riemann data leading to nucleation or cavitation, i.e. the creation of a new phase.
We conclude this work with Section \ref{chap_conclusion} where we give a detailed discussion of the assumptions made to state the previous results followed by some examples and the conclusion.

    \newpage
    \section{Isothermal Euler Equations}\label{chap_euler_eq}
\setcounter{para}{1}
In this work we study inviscid, compressible and isothermal two phase flows.
The two phases are either the liquid or the vapor phase of one substance. The phases are distinguished by the \textit{mass density} $\rho$ and further described by the \textit{velocity} $u$.
Sometimes it is convenient to use the \textit{specific volume} $v = 1/\rho$ instead of the mass density. We will make the reader aware of such situations.
The physical quantities depend on time $t \in \mathbb{R}_{\geq 0}$ and space $x \in \mathbb{R}$.
In regular points of the bulk phases we have the conservation law for mass and the balance law for momentum, i.e.
\begin{align}
    \partial_t\rho + \partial_x(\rho u) &= 0,\label{mass_cons}\\
    \partial_t(\rho u) + \partial_x(\rho u^2 + p) &= 0.\label{mom_balance}
\end{align}
The system of equations (\ref{mass_cons}) and (\ref{mom_balance}) is referred to as the \textit{isothermal Euler equations}.
The additional quantity $p$ denotes the \textit{pressure} and is related to the mass density via the \textit{equation of state (EOS)} $p = p(\rho)$.
Sometimes one also refers to the EOS as \textit{pressure law}. Such an EOS crucially depends on the considered substance and how this substance is modeled.
Across any discontinuity we have the following jump conditions
\begin{align}
    [\![\rho(u - W)]\!] &= 0,\label{jc_mass}\\
    \rho(u - W)[\![u]\!] + [\![p]\!] &= 0.\label{jc_mom}
\end{align}
Here we write $[\![\Psi]\!] = \Psi^+ - \Psi^-$, where $\Psi^+$ is the right and $\Psi^-$ the left sided limit of the physical quantity $\Psi$.
Furthermore every discontinuity satisfies the following \textit{entropy inequality}
\begin{align}
    \rho(u - W)[\![g + e^{kin}]\!] \leq 0.\label{entropy_ineq}
\end{align}
Further, $W$ denotes the speed of the discontinuity and $Z = -\rho(u - W)$ the \textit{mass flux} where we will distinguish between a classical shock wave and the phase boundary (non-classical shock)
\begin{align*}
    Z = \begin{cases} Q,\quad\text{shock wave}\\ z,\quad\text{phase boundary}\end{cases}
    \quad\text{and}\quad
    W = \begin{cases} S,\quad\text{shock wave}\\ w,\quad\text{phase boundary}\end{cases}.
\end{align*}
\subsection{Definition and Requirements for the \textit{EOS}}\label{sec_eos}
Usually one only works with the pressure law when dealing with the Euler equations. Nevertheless the pressure law does not contain all the information about a fluid
or more general a thermodynamic system.
From a thermodynamic point of view a system in (local) equilibrium can be described relating the extensive quantities \textit{energy} $E$, \textit{volume} $V$ and \textit{entropy} $S$, i.e. $E(V,S)$.
In the following we will use the corresponding (intensive) densities and thus we use small letters (e.g.: $e,v,s$).
Given this relation every other quantity can be derived using the first and second law of thermodynamics and the so called Maxwell relations.
A condensed overview, including the difference between a complete and an incomplete EOS, can be found in \cite{Menikoff1989}.
For detailed information about EOS we refer to standard literature, cf. \cite{Bartelmann2015, Landau1987, Mueller1985, Mueller2009, Rebhan2010}.
A discussion using the ideal gas EOS and the Tait EOS can also be found in \cite{Dreyer2012}. From this point on we assume that we have an EOS for each phase with consistent thermodynamic properties.
There are different possible thermodynamic potentials which can be used to describe a system and they are all connected to each other using the Legendre transform.
Thus one can start from any potential and will get similar results.
For the discussion of the equations at the interface we need the Gibbs energy and hence shortly summarize the most important features, i.e.\ those we need for our purpose.
More details can be found in the above mentioned literature and references therein.
\begin{defi}[Gibbs Energy and Sound Speed]\label{defi_gibbs}
    The Gibbs energy is a function of the pressure $p$ and the temperature $T$. The (complete) differential is given by
    \begin{align*}
        \dd g = -s\dd T + v\dd p.
    \end{align*}
    Further we define the isothermal sound speed as
    \begin{align*}
        a = \sqrt{-v^2\left(\dfrac{\partial p}{\partial v}\right)_T}.
    \end{align*}
\end{defi}
From Definition \ref{defi_gibbs} we obtain
\begin{align}
    \left(\dfrac{\partial g}{\partial p}\right)_T &= v > 0,\quad \left(\dfrac{\partial^2 g}{\partial p^2}\right)_T = \left(\dfrac{\partial v}{\partial p}\right)_T = -\left(\frac{v}{a}\right)^2 < 0.
    \label{eq:gibbs_derivatives}
\end{align}
Since thermodynamic quantities may be expressed using different choices of independent variables the brackets with the subscript simply denote which quantity is held constant when calculating the
derivative. In the isothermal case the Gibbs potential just depends on the pressure and hence we omit writing the brackets with subscript $T$.
Here the volume $v$ and the speed of sound $a$ are strictly positive functions of the pressure $p$.
Furthermore the inequality for the second derivative is due to the requirement of thermodynamic stability for an isothermal system.
In short this can be seen by considering the requirements for the full case.
There, thermodynamic stability requires the energy to be a convex function, both in the entropy and the volume. This implies that the Hessian of the energy is non negative.
If we now assume the temperature to be constant, what remains is
\begin{align}
    0 \leq \frac{\dd^2 e}{\dd v^2} = -\frac{\dd p}{\dd v}.\label{thermodyn_stab}
\end{align}
%
In the following we use the subscripts $\{V,L\}$ when it is necessary to distinguish the vapor and the liquid phase.
Since we are concerned with two phases we write $g_L$ for the Gibbs energy of the liquid phase and $g_V$ for the vapor phase, respectively. Further we require
\begin{align*}
    \dfrac{\partial g_i}{\partial p_j} &= 0,\,\, i \neq j,\,\, i, j \in \{V, L\}.
\end{align*}
Since we only consider one substance the condition for two phases to be in equilibrium is
\begin{align}
    g_L(p_L) = g_V(p_V).\label{defi_equil}
\end{align}
Due to the monotonicity of $g_K,\,\, K \in \{V,L\}$ we have
\begin{align*}
    g_L(p_L) = g_V(p_V)\quad\Leftrightarrow\quad p_L = p_V
\end{align*}
and we write in this case
\begin{align*}
    p_L = p_V =: p_0\quad\text{and}\quad g_L(p_0) = g_V(p_0).
\end{align*}
A crucial point when dealing with different phases is how to discriminate them and how to connect them thermodynamically consistent.
Equations of state describing two phases (e.g. \textit{van der Waals} EOS) have a so called spinodal region which is avoided by the Maxwell construction (or equal area rule).
We want to discriminate the phases using the specific volumes. Therefore we need an upper bound for the liquid volume $v_m$ and a lower bound for the vapor volume $\tilde{v}$.
This should still be consistent with the Maxwell construction.
Therefore we may proceed as follows. We use the EOS for each phase and prescribe the minimum liquid pressure $p_{min}$ (e.g. $p_{min} = 0$) and from this we obtain $v_m$.
Further we know the saturation pressure $p_0$ for a given temperature $T_0$ from a calculation or from tables which are available for many substances, such as for water \cite{Wagner1998}.
Now we connect our two EOS monotonically and then obtain the maximum vapor pressure $\tilde{p}$ using the Maxwell construction, see \cite{Mueller2009}.
\begin{defi}[Maximum Vapor Pressure]\label{max_pV}
    Given a fixed temperature $T_0$ the corresponding saturation pressure $p_0$ is given by (\ref{defi_equil}). Furthermore $p_{min}$ is defined to be the minimum liquid pressure.
    Let $\bar{v}(p)$ be a function such that
    \begin{align*}
        v_L(p_{min}) = \bar{v}(p_{min}),\quad v_V(\tilde{p}) = \bar{v}(\tilde{p})\quad\text{and}\quad \bar{v}'(p) > 0.
    \end{align*}
    Then the \textit{maximum vapor pressure} $\tilde{p}$ is found as the solution of the following equation
    \begin{align*}
        0 = p_0(v_V(p_0) - v_L(p_0)) - \int_{v_L(p_0)}^{v_V(p_0)} p(v)\,\textup{d}v.
    \end{align*}
    The function $p(v)$ given by
    \begin{align*}
        p(v) = \begin{cases} p_L(v),\,v \in \bigl(0,v_L(p_{min})\bigr]\\ \bar{p}(v),\,v \in \bigl(v_L(p_{min}),v_V(\tilde{p})\bigr),\\ p_V(v),\,v \in \bigl[v_V(\tilde{p}),\infty\bigr) \end{cases}.
    \end{align*}
\end{defi}
Finally, analogous to \cite{Menikoff1989} we introduce dimensionless quantities which we will use later on.
\begin{defi}[Dimensionless Quantities]\label{def_gamma_G}
    We define the (isothermal) dimensionless speed of sound as
    \begin{align*}
        \gamma := -\frac{v}{p}\frac{\dd p}{\dd v}.
    \end{align*}
    and the (isothermal) fundamental derivative
    \begin{align*}
        \mathcal{G} := -\frac{1}{2}v\dfrac{\dfrac{\dd^2p}{\dd v^2}}{\dfrac{\dd p}{\dd v}}.
    \end{align*}
\end{defi}
It is straight forward to verify and no surprise that these quantities are completely analogue to those defined in \cite{Menikoff1989}.
In fact, by using the relations 
given in \cite{Menikoff1989} and assuming the temperature to be fixed, one also obtains the results given above.
However we want to emphasize that $\gamma$ and $\mathcal{G}$ defined here are \textit{not equal} to those defined in \cite{Menikoff1989}.
This is because we assume the temperature to be constant, whereas in \cite{Menikoff1989} the derivatives are taken at constant entropy.
To clarify this, let us for the moment write 
$\gamma_S$ for the isentropic quantity defined in \cite{Menikoff1989}. Then we have (cf. \cite{Menikoff1989})
\begin{align*}
    \frac{\gamma}{\gamma_S} = \frac{c_V}{c_p}
\end{align*}
and hence $\gamma \leq \gamma_S$ for thermodynamic stable systems. 
Further we have for $\gamma$, using Definition \ref{defi_gibbs}
\begin{align}
    \gamma = \frac{a^2}{pv}.\label{dimensionless_a}
\end{align}
For the fundamental derivative one may also write
\begin{align}
    \mathcal{G} = \frac{1}{2}\frac{v^2}{p\gamma}\frac{\dd^2p}{\dd v^2} = -\frac{v}{a}\frac{\dd a}{\dd v} + 1\label{fund_deriv_1}
\end{align}
or when expressed in terms of the pressure
\begin{align}
    \mathcal{G} = \frac{a}{v}\frac{\dd a}{\dd p} + 1.\label{fund_deriv_2}
\end{align}
The \textit{isotherms} in the $p-v$ plane are convex if $\mathcal{G} > 0$, which we will assume from now on.
\subsection{Riemann Problem}\label{sec_riemann_problem}
In the following we briefly discuss the solution of the Riemann problem for the isothermal Euler equations (\ref{mass_cons})-(\ref{mom_balance}) for a single phase.
In order to do so we will discuss the elementary wave types that can occur, which are shock or rarefaction waves.
The Riemann problem is given by equations (\ref{mass_cons})-(\ref{mom_balance}), the EOS and the Riemann initial data
\begin{align}
    \rho(x,0) = \begin{cases} &\rho_-,\,x < 0\\ &\rho_+,\,x > 0\end{cases}
    \quad\text{and}\quad
    u(x,0) = \begin{cases} &u_-,\,x < 0\\ &u_+,\,x > 0\end{cases}.\label{init_data_rp_0}
\end{align}
We rewrite the system (\ref{mass_cons}) - (\ref{mom_balance}) in quasilinear form in terms
of the primitive variables, i.e. the density $\rho$ and the velocity $u$
\begin{align}
    \left(\begin{matrix} \rho\\ u \end{matrix}\right)_t
    +
    \left(\begin{matrix}
        u                 & \rho\\
        \frac{a^2}{\rho}  & u
    \end{matrix}\right)
    \left(\begin{matrix} \rho\\ u \end{matrix}\right)_x
    = 0.\label{eulereq_quasilin}
\end{align}
The Jacobian matrix
\begin{align}
    \mathbf{A} =
    \left(\begin{matrix}
        u                 & \rho\\
        \frac{a^2}{\rho}  & u
    \end{matrix}\right)\label{eulereq_jacobian}
\end{align}
has the following eigenvalues and corresponding eigenvectors
\begin{align}
    \lambda_1 = u - a,\quad\mathbf{r}_1 = \left(\begin{matrix} \rho\\ -a \end{matrix}\right),
    \quad
    \lambda_2 = u + a,\quad\mathbf{r}_2 = \left(\begin{matrix} \rho\\ a \end{matrix}\right).\label{jacobian_eigen_val_vect}
\end{align}
Due to the requirement of thermodynamic stability (\ref{thermodyn_stab}) this system is hyperbolic. We have strict hyperbolicity for
\begin{align}
    \gamma > 0.\label{strict_thermodyn_stab}
\end{align}
Furthermore one can immediately verify that the waves corresponding to the eigenvalues and
eigenvectors are genuine nonlinear if and only if the fundamental derivative
\begin{align*}
    \mathcal{G} &= \frac{\rho}{a}\frac{\dd a}{\dd \rho} + 1.
\end{align*}
does not vanish, i.e.
\begin{align}
    \nabla \lambda_{1/2}\cdot\mathbf{r}_{1/2} = \mp\frac{a}{\rho}\mathcal{G} \neq 0.
\end{align}
Here this is in fact the case, since we assumed $\mathcal{G} > 0$.
For systems with genuine nonlinear waves the Lax condition is enough to pick the right solution, cf. \cite{LeFloch2002} and also \cite{Menikoff1989} for the full system.
The Riemann invariants for this system are
\begin{align}
    I_1 = u + \int\frac{a}{\rho}\,\dd\rho\quad\text{and}\quad I_2 = u - \int\frac{a}{\rho}\,\dd\rho.\label{riem_inv}
\end{align}
\subsubsection{Entropy Inequality across a Shock Wave}\label{sec_entropy_ineq}
Hantke et al.\ proved, that the Lax condition is equivalent to the entropy condition for an isothermal system.
This holds true for the general entropy inequality given by (\ref{entropy_ineq})
\begin{align*}
    Q[\![g + e^{kin}]\!] = -\rho(u - S)[\![g + e^{kin}]\!] \geq 0.
\end{align*}
Consider two states
\begin{align*}
    \left(\begin{matrix} \rho_1\\ u_1 \end{matrix}\right)\quad\text{and}\quad\left(\begin{matrix} \rho_2\\ u_2 \end{matrix}\right)
\end{align*}
separated by a shock wave moving with speed $S$. Using the specific volume $v = 1/\rho$ one obtains
\begin{align}
    \frac{a(p_1)^2}{v(p_1)^2} < Q^2 < \frac{a(p_2)^2}{v(p_2)^2}.\label{Q_square_ineq}
\end{align}
which gives the Lax condition for a left Shock ($Q > 0$) and a right shock ($Q < 0$).
\subsubsection{Rarefaction Wave}\label{sec_rarefaction_wave}
For a rarefaction wave we use the Riemann invariants (\ref{riem_inv}) and hence obtain for a left rarefaction wave (corresponding to $\lambda_1$)
\begin{align}
    u_2 - u_1 = -\int_{\rho_1}^{\rho_2} \frac{a}{\rho}\,\dd\rho.\label{left_raref}
\end{align}
Furthermore the slope inside the rarefaction is given by
\begin{align}
    \frac{\dd x}{\dd t} = \frac{x}{t} = \lambda_1 = u - a\label{left_raref_slope}
\end{align}
and hence we obtain for the solution inside the rarefaction fan
\begin{align}
    u = \frac{x}{t} + a\quad\text{and}\quad F(\rho) = u - u_1 + \int_{\rho_1}^{\rho} \frac{a}{\sigma}\,\dd\sigma = 0.\label{left_raref_fan}
\end{align}
Here $\rho$ is obtained as the root of $F(\rho)$. Similar we obtain the results for a right rarefaction
\begin{align}
    u_2 - u_1 &= \int_{\rho_1}^{\rho_2} \frac{a}{\rho}\,\dd\rho,\quad\frac{\dd x}{\dd t} = \frac{x}{t} = \lambda_2 = u + a,\notag\\
    u &= \frac{x}{t} - a\quad\text{and}\quad F(\rho) = u_2 - u - \int_{\rho}^{\rho_2} \frac{a}{\sigma}\,\dd\sigma = 0.\label{right_raref}
\end{align}
\subsubsection{Shock Wave}\label{sec_shock_wave}
The relation across a shock wave is given by 
\begin{align}
    [\![u]\!]^2 = -\llbracket p \rrbracket[\![v]\!] = \frac{[\![p]\!][\![\rho]\!]}{\rho_1\rho_2}
    \quad\Leftrightarrow\quad
    [\![u]\!] = -\sqrt{-\llbracket p \rrbracket[\![v]\!]} = -\sqrt{\frac{[\![p]\!][\![\rho]\!]}{\rho_1\rho_2}}.\label{shock_function_gen}
\end{align}
\subsubsection{Solution of the Riemann Problem}\label{sec_rp_solution}
If we now want to solve the Riemann problem for the isothermal Euler equations we just have to connect the three constant states separated by the waves using the equations obtained above.
Therefore we assume the left and right state to be given and use that the velocity between the waves is constant. The solution is obtained as the root of the following function
\begin{align}
    f(\rho,W_L,W_R) &= f_R(\rho,W_R) + f_L(\rho,W_L) + u_R - u_L = 0,\label{1phase_rp_sol_rho}\\
    f_K(\rho,W_K) &=
    \begin{dcases}
        \sqrt{\frac{[\![p]\!][\![\rho]\!]}{\rho\rho_K}},\,\,\rho > \rho_K\,\,\text{(Shock)}\\
        \int_{\rho_K}^\rho \dfrac{a(\sigma)}{\sigma}\,\dd\sigma,\,\,\rho \leq \rho_K\,\,\text{(Rarefaction)}
    \end{dcases},\quad K \in \{L, R\}.\notag
\end{align}
Due to $p'(\rho) > 0$ we could also state this problems in terms of the unknown pressure $p$, i.e.
\begin{align}
    f(p,W_L,W_R) &= f_R(p,W_R) + f_L(p,W_L) + u_R - u_L = 0,\label{1phase_rp_sol_p}\\
    f_K(p,W_K) &=
    \begin{dcases}
        \sqrt{-[\![p]\!][\![v]\!]},\,\,p > p_K\,\,\text{(Shock)}\\
        \int_{p_K}^p \frac{v(\zeta)}{a(\zeta)}\,\textup{d}\zeta,\,\,p \leq p_K\,\,\text{(Rarefaction)}
    \end{dcases},\quad K \in \{L, R\}.\notag
\end{align}
In order to investigate $f(p,W_L,W_R)$ we need information about the asymptotic behavior 
\begin{align*}
    v(p) \stackrel{p\to\infty}\to 0,\quad v(p) \stackrel{p \to 0}{\to} \infty
    \quad\text{and further}\quad
    \dfrac{\dd v(p)}{\dd p} \stackrel{(\ref{eq:gibbs_derivatives})_2}{=} -\frac{v(p)^2}{a(p)^2} < 0.
\end{align*}
We obtain for $f_K(p,W_K)$ in the case of a shock wave
\begin{align}
    \frac{\dd}{\dd p}f_K(p,W_K) &= \frac{-[\![v]\!] + [\![p]\!]\frac{v^2}{a^2}}{2\sqrt{-[\![p]\!][\![v]\!]}} > 0,\label{mono_conc_shock}\\
    \frac{\dd^2}{\dd p^2}f_K(p,W_K)
    &= -\frac{1}{4(-[\![p]\!][\![v]\!])^{3/2}}\left(-4[\![p]\!]^2[\![v]\!]\frac{v^3}{a^4}\mathcal{G} + \left([\![p]\!]\frac{v^2}{a^2} - [\![v]\!]\right)^2\right) < 0\notag
\end{align}
For a rarefaction wave we yield
\begin{align}
    \frac{\dd}{\dd p}f_K(p,W_K) &= \frac{v(p)}{a(p)} > 0,\label{mono_conc_raref}\\
    \frac{\dd^2}{\dd p^2}f_K(p,W_K) &= -\frac{v(p)^2}{a(p)^3}\mathcal{G} < 0\notag
\end{align}
Combining (\ref{mono_conc_shock}) with (\ref{mono_conc_raref}) gives
\begin{align}
    \frac{\dd}{\dd p}f(p,W_L,W_R) > 0\quad\text{and}\quad\frac{d^2}{dp^2}f(p,W_L,W_R) < 0.\label{mono_conc_rp_fun}
\end{align}
Using the asymptotic behavior of $v(p)$ gives
\begin{align}
    f(p,W_L,W_R) \stackrel{p \to 0}{\to} -\infty\quad\text{and}\quad f(p,W_L,W_R) \stackrel{p \to \infty}{\to} +\infty\label{asymp_rp_fun}
\end{align}
and hence we have a unique root which determines the solution of our system.

    \newpage
    \section{Solution at the Interface}\label{chap_pb_sol}
\setcounter{para}{1}
The phase boundary separating the liquid and the vapor phase is a non-classical or under compressive shock, see \cite{Dafermos2010} or \cite{LeFloch2002} and references therein.
Hence the Lax criterion alone will not give us a unique solution and we need a further relation at the interface. This equation is called \textit{kinetic relation}.
We use the kinetic relation derived by Dreyer et al. \cite{Dreyer2012}.
The kinetic relation is chosen such that the mass flux $z$ is proportional to the jump term in the entropy inequality (\ref{entropy_ineq})
\begin{align*}
    z[\![g + e^{kin}]\!] \geq 0.
\end{align*}
If we assume the vapor left to the liquid phase the kinetic relation reads
\begin{align}
    z = \tau p_V[\![g + e^{kin}]\!] = \tau p_V[g_L - g_V + e^{kin}_L - e^{kin}_V].\label{kin_rel1}
\end{align}
Otherwise we can use
\begin{align}
    z = -\tau p_V[g_L - g_V + e^{kin}_L - e^{kin}_V].\label{kin_rel_lv}
\end{align}
In the following we will assume the first case.
In this section we will prove that there exists a unique solution of the equations at the interface provided certain conditions hold.
By this we mean that there exists a unique liquid (vapor) state for a prescribed vapor (liquid) state such that the following equations hold
\begin{align*}
    [\![z]\!] &= 0,\\
    -z[\![u]\!] + [\![p]\!] &= 0,\\
    z &= \tau p_V[\![g + e^{kin}]\!].
\end{align*}
Here $e^{kin}$ denotes the \textit{kinetic energy}. Furthermore we have for the so called \textit{mobility} $0 < \tau \in \mathbb{R}$. Usually one uses
\begin{align}
    \tau = \frac{1}{\sqrt{2\pi}}\left(\frac{m}{kT_0}\right)^{\frac{3}{2}}\label{tau_special}
\end{align}
where $m$ denotes the mass of a single molecule, $k$ the Boltzmann constant and $T_0$ the fixed temperature, see \cite{Bond2004,Dreyer2012}.
Using the jump conditions (\ref{jc_mass})-(\ref{jc_mom}) we can rewrite (\ref{kin_rel1}) and obtain
\begin{align}
    z = \tau p_V[\![g - \frac{1}{2}p(v_L + v_V)]\!].\label{kin_rel2}
\end{align}
Furthermore we can combine the jump conditions and obtain
\begin{align}
    [\![p]\!] + z^2[\![v]\!] = 0.\label{pb_mom_balance}
\end{align}
Together with the EOS and (\ref{kin_rel2}) equation (\ref{pb_mom_balance}) is a single equation for one unknown given one state at the phase boundary.
For example we will prescribe the vapor pressure and then obtain the liquid pressure as the solution of equation (\ref{pb_mom_balance}).
In the following we will assume as before that $\gamma_V \geq 0$ and $\mathcal{G}_K > 0,\,K \in \{V,L\}$.
From the mathematical point of view we need further assumptions to solve the problem.
A discussion will be given later on and it will turn out that these assumptions are rather liberal from a physical point of view, see Subsection \ref{discuss_assumptions}.
In the following we need the quotient of the specific volumes to be uniformly bounded as well as the corresponding sound speeds
\begin{align}
    0 &< \frac{v_L}{v_V} \leq \alpha < 1,\quad 0 < \frac{v_L}{v_V}\frac{a_V}{a_L} \leq \alpha\beta < 1,\quad \tau(1 - \alpha)^2a_V^3 < \gamma_V\quad\text{and}\notag\\
    0 &< p_V \leq \sigma_{max}p_0\quad\text{with}\quad\sigma_{max} = \frac{1 + \sqrt{11 - 6\alpha}}{2}.\label{assumptions1}
\end{align}
\begin{rem}\label{rem:p_dependence}
    The specific volume and the speed of sound depend on the pressure but for convenience we often will not write out this dependence explicitly.
\end{rem}
Now we can state one of the main results of this work.

\begin{thm}[Existence and Uniqueness of a Solution at the Interface]\label{exis_uniq_thm}
    For two phases each described by a thermodynamic consistent equation of state meeting the requirements (\ref{assumptions1}) and
    \begin{align*}
        -a_V/v_V \leq z \leq a_L/v_L
    \end{align*}
    exists a unique solution of equation (\ref{pb_mom_balance}).
    Furthermore the mass flux $z$ is uniquely defined. The liquid pressure can be written as a function of the vapor pressure and has the following properties
    \begin{align*}
        p_L^\ast = \varphi(p_V^\ast) \geq p_V^\ast,\,\,\varphi(p_0) = p_0,\,\,\frac{\textup{d}\varphi(p_V^\ast)}{\textup{d}p_V^\ast} > 0
    \end{align*}
\end{thm}
In the remaining part of this section we will give the proof of this theorem.
\subsection{Proof}\label{sec_exis_uniq_proof}
The proof of Theorem \ref{exis_uniq_thm} is based on the \textit{implicit function theorem}. The main steps are the following
\begin{enumerate}[(i)]
    \item We define a function $f(p_V,p_L)$, see (\ref{math_fun_mom_balance}), which we will analyze and where the roots correspond to the solution of (\ref{pb_mom_balance}).
    \item The local existence of an admissible root, see Definition \ref{adm_sol}, for the equilibrium case $(p_0,p_0)$ is given in Remark \ref{rem:equil_solution}.
    \item Lemma \ref{f_pL_pos} and Lemma \ref{f_pV_neg} state that the first order derivatives of $f(p_V,p_L)$ each have a sign for an admissible solution.
    \item Uniqueness is shown in Lemma \ref{unique_sol} and global existence is stated and proven in Lemma \ref{glob_exist}.
\end{enumerate}
Replacing $z$ in (\ref{pb_mom_balance}) using (\ref{kin_rel2}) we obtain
\begin{align*}
    [\![p]\!] + \left(\tau p_V[\![g - \frac{1}{2}p(v_L + v_V)]\!]\right)^2[\![v]\!] = 0.
\end{align*}
According to this equation we define the following functions
\begin{align}
    h(p_V,p_L) &:= \tau[\![g - \frac{1}{2}p(v_L + v_V)]\!]\notag\\
               &= \tau\left[g_L(p_L) - g_V(p_V) - \frac{1}{2}(p_L - p_V)(v_L(p_L) + v_V(p_V))\right],\notag\\
    f(p_V,p_L) &:= [\![p]\!] + \left(p_Vh(p_V,p_L)\right)^2[\![v]\!].\label{math_fun_mom_balance}
\end{align}
Obviously every root of (\ref{math_fun_mom_balance}) is a solution of (\ref{pb_mom_balance}) and we easily see
\begin{align}
    0 = f(p_V^\ast,p^\ast_L) \quad\stackrel{[\![v]\!] < 0}{\Rightarrow}\quad [\![p]\!] \geq 0.\label{p_jump_pos}
\end{align}
Let us furthermore define the following
\begin{defi}[Admissible Solution]\label{adm_sol}
    Let $(p_V^\ast,p_L^\ast)$ be a solution of $f(p_V^\ast,p_L^\ast) = 0$. We say this solution is \textit{admissible} if further the following inequalities hold
    \begin{align*}
        -\frac{a_V(p_V^\ast)}{v_V(p_V^\ast)} \leq p_V^\ast h(p_V^\ast,p_L^\ast) \leq \frac{a_L(p_L^\ast)}{v_L(p_L^\ast)}.
    \end{align*}
\end{defi}
The quantities $a_K$ and $v_K$ with $K \in \{L,V\}$ are functions of the pressure as already mentioned in Remark \ref{rem:p_dependence}. Thus the bounds are evaluated at the pressures
$(p_V^\ast,p_L^\ast)$ which solve $f(p_V^\ast,p_L^\ast) = 0$.
\begin{rem}\label{rem:equil_solution}
    It is immediately verified that a solution $f(p_V^\ast,p_L^\ast) = 0$ with $p_V^\ast = p_L^\ast =: p_0$ implies equilibrium $g_L(p_L^\ast) = g_V(p_V^\ast)$ and vice versa.
    Thus we further obtain
    \begin{align}
        f(p_0,p_0) = 0,\,\,\partial_{p_V}f(p_0,p_0) = -1,\,\,\partial_{p_L}f(p_0,p_0) = 1\quad\text{with}\quad p_0h(p_0,p_0) = 0.\label{equil_sol}
    \end{align}
    Hence there exists a neighborhood of $p_V = p_0$ such that (\ref{pb_mom_balance}) implicitly defines a function $p_L = \varphi(p_V)$ with
    $\varphi'(p_V) > 0$. Additionally $(p_0, p_0)$ is an admissible solution with $z = 0$.
\end{rem}
\begin{lem}\label{h_pL_decrease}
    The function $h(p_V,p_L)$ is strictly monotonically decreasing in $p_L$ under the given assumptions, i.e.
    \begin{align*}
        \partial_{p_L}h(p_V,p_L) < 0.
    \end{align*}
\end{lem}
\begin{prf}
    We obtain for the partial derivative of $h(p_V,p_L)$ using $(\ref{eq:gibbs_derivatives})_2$
    \begin{align*}
        \partial_{p_L}h(p_V,p_L) = \frac{\tau}{2}\left\{[\![v]\!] + [\![p]\!]\frac{v_L^2}{a_L^2}\right\}.
    \end{align*}
    Let us consider $[\![p]\!] \geq 0$ since it is the only relevant case and the statement is obvious for $[\![p]\!] \leq $ anyway.
    Since $\mathcal{G}_L > 0$ we yield for the second partial derivative with respect to $p_L$ using $(\ref{eq:gibbs_derivatives})_2$ and (\ref{fund_deriv_2})
    \begin{align*}
        \partial^2_{p_L}h(p_V,p_L) = -\tau[\![p]\!]\frac{v_L^3}{a_L^4}\mathcal{G}_L < 0.
    \end{align*}
    For $p_L = p_V$ we know that the Lemma is true and if we increase $p_L$ the function is decreasing. Keep in mind that we have $[\![p]\!] > 0$.
    Hence we conclude $\partial_{p_L}h(p_V,p_L) < 0$.
\end{prf}
\begin{cor}\label{z_simple_upbound}
    Every root of (\ref{math_fun_mom_balance}) with $z > 0$ is admissible.
\end{cor}
\begin{prf}
    Using Lemma \ref{h_pL_decrease} one obtains for $f(p_V^\ast,p_L^\ast) = 0$ with $z = p^\ast_Vh(p_V^\ast,p^\ast_L)$
    \begin{align*}
        z^2 = \left(p^\ast_Vh(p_V^\ast,p^\ast_L)\right)^2 \stackrel{(\ref{math_fun_mom_balance})}{=} -\frac{[\![p^\ast]\!]}{[\![v(p^\ast)]\!]}
        \stackrel{\text{Lemma} \ref{h_pL_decrease}}{<} \frac{a_L^2}{v_L^2}.
    \end{align*}
\end{prf}
\begin{lem}\label{f_pL_pos}
    Let $(p_V^\ast,p_L^\ast)$ be an admissible solution of $f(p_V^\ast,p_L^\ast) = 0$. Then the following inequality holds
    \begin{align*}
        \partial_{p_L}f(p_V^\ast,p_L^\ast) > 0.
    \end{align*}
\end{lem}
\begin{prf}
    For the equilibrium solution (\ref{equil_sol}) the stated relation is obvious. Let us consider $p_V^\ast h(p_V^\ast,p_L^\ast) > 0$.
    Using Lemma \ref{h_pL_decrease} and $[\![v]\!] < 0$ we have
    \begin{align*}
        \partial_{p_L}f(p_V^\ast,p_L^\ast) = 1 + 2\underbrace{(p_V^\ast h(p_V^\ast,p_L^\ast))}_{> 0}\underbrace{(p_V^\ast\partial_{p_L}h(p_V^\ast,p_L^\ast))[\![v]\!]}_{> 0}
                                           - \underbrace{(p_V^\ast h(p_V^\ast,p_L^\ast))^2\frac{v_L^2}{a_L^2}}_{< 1} > 0.
    \end{align*}
    It remains to prove the Lemma for the case $p_V^\ast h(p_V^\ast,p_L^\ast) < 0$. We can write
    \begin{align*}
        \partial_{p_L}f(p_V^\ast,p_L^\ast) &= 1 + 2(p_V^\ast h(p_V^\ast,p_L^\ast))(p_V^\ast\partial_{p_L}h(p_V^\ast,p_L^\ast))[\![v]\!] - (p_V^\ast h(p_V^\ast,p_L^\ast))^2\frac{v_L^2}{a_L^2}\\
        &= 1 + \tau p_V^\ast(p_V^\ast h(p_V^\ast,p_L^\ast))[\![v]\!]^2\left(1 - (p_V^\ast h(p_V^\ast,p_L^\ast))^2\frac{v_L^2}{a_L^2}\right) - (p_V^\ast h(p_V^\ast,p_L^\ast))^2\frac{v_L^2}{a_L^2}\\
        &= \left(1 - (p_V^\ast h(p_V^\ast,p_L^\ast))^2\frac{v_L^2}{a_L^2}\right)\left(1 + \tau p_V^\ast(p_V^\ast h(p_V^\ast,p_L^\ast))[\![v]\!]^2\right).
    \end{align*}
    The first term is positive, because of $-a_V/v_V \leq p_V^\ast h(p_V^\ast,p_L^\ast) < 0$ and $a_V^2/v_V^2 < a_L^2/v_L^2$. For the second term we have
    \begin{align*}
        &0 < 1 + \tau p_V^\ast(p_V^\ast h(p_V^\ast,p_L^\ast))[\![v]\!]^2\quad\stackrel{p_V^\ast h(p_V^\ast,p_L^\ast) < 0}{\Leftrightarrow}\quad
        \tau < -\frac{1}{p_V^\ast(p_V^\ast h(p_V^\ast,p_L^\ast))[\![v]\!]^2}.
    \end{align*}
    Indeed we obtain
    \begin{align*}
        -&\frac{1}{p_V^\ast(p_V^\ast h(p_V^\ast,p_L^\ast))[\![v]\!]^2} > \frac{v_V}{p_V^\ast a_V[\![v]\!]^2} = \frac{1}{p_V^\ast v_V a_V \left(\dfrac{v_L}{v_V} - 1\right)^2}
        \stackrel{(\ref{assumptions1})_1}{\geq} \frac{\gamma_V}{(1 - \alpha)^2a_V^3} \stackrel{(\ref{assumptions1})_3}{>} \tau.
    \end{align*}
    This proves the Lemma.
\end{prf}
\begin{lem}\label{f_pV_neg}
    Let $(p_V^\ast,p_L^\ast)$ be an admissible solution of $f(p_V^\ast,p_L^\ast) = 0$. Then the following inequality holds
    \begin{align}
        \partial_{p_V}f(p_V^\ast,p_L^\ast) < 0.
    \end{align}
\end{lem}
\begin{prf}
    Since we have $f(p_V^\ast,p_L^\ast) = 0$ we can write for $\partial_{p_V}h(p_V^\ast,p_L^\ast)$
    \begin{align}
        \partial_{p_V}h(p_V^\ast,p_L^\ast) = \frac{\tau}{2}\left\{[\![v]\!] + [\![p]\!]\frac{v_V^2}{a_V^2}\right\}
                                           \stackrel{(\ref{math_fun_mom_balance})}{=} \frac{\tau}{2}[\![v]\!]\left(1 - (p_V^\ast h(p_V^\ast,p_L^\ast))^2\frac{v_V^2}{a_V^2}\right)\label{dpvh}
    \end{align}
    and hence we conclude
    \begin{align*}
        \partial_{p_V}h(p_V^\ast,p_L^\ast)\,
        \begin{cases}
            &< 0,\quad (p_V^\ast h(p_V^\ast,p_L^\ast))^2 < \dfrac{a_V^2}{v_V^2},\\
            &\geq 0,\quad (p_V^\ast h(p_V^\ast,p_L^\ast))^2 \geq \dfrac{a_V^2}{v_V^2}.
        \end{cases}
    \end{align*}
    In the following we will discuss three cases depending on $p_V^\ast h(p_V^\ast,p_L^\ast)$.\\
    \newline
    \textbf{First Case:} We discuss the case where $-a_V/v_V \leq p_V^\ast h(p_V^\ast,p_L^\ast) \leq 0$. It is obvious to see
    \begin{align*}
        \partial_{p_V}f(p_V^\ast,p_L^\ast) = 
        \left.\begin{cases}
            -1&,\quad p_V^\ast h(p_V^\ast,p_L^\ast) = 0,\\
            2\dfrac{a_V^2}{p_V^\ast v_V^2}[\![v]\!]&,\quad p_V^\ast h(p_V^\ast,p_L^\ast) = -\dfrac{a_V}{v_V}
        \end{cases}\right\} < 0.
    \end{align*}
    In between we have $-a_V/v_V < p_V^\ast h(p_V^\ast,p_L^\ast) < 0$ and so all together
    \begin{align*}
        &\partial_{p_V}f(p_V^\ast,p_L^\ast) =\dots\\
        &= -1 + 2\underbrace{(p_V^\ast h(p_V^\ast,p_L^\ast))}_{< 0}\underbrace{(h(p_V^\ast,p_L^\ast) + p_V^\ast\partial_{p_V}h(p_V^\ast,p_L^\ast))}_{< 0}[\![v]\!] 
         + \underbrace{(p_V^\ast h(p_V^\ast,p_L^\ast))^2\frac{v_V^2}{a_v^2}}_{< 1} < 0.
    \end{align*}
    For $0 < p_V^\ast h(p_V^\ast,p_L^\ast) < a_L/v_L$ we split the proof into two parts. First we discuss the interval up to $a_V/v_V$ and then the remaining part smaller than $a_L/v_L$.\\
    \newline
    \textbf{Second Case:} Using $0 < p_V^\ast h(p_V^\ast,p_L^\ast) \leq a_V/v_V$ we obtain
    \begin{align*}
        &\partial_{p_V}f(p_V^\ast,p_L^\ast) =\dots\\
        &= -1 + 2(p_V^\ast h(p_V^\ast,p_L^\ast))(h(p_V^\ast,p_L^\ast) + p_V^\ast\partial_{p_V}h(p_V^\ast,p_L^\ast))[\![v]\!] + (p_V^\ast h(p_V^\ast,p_L^\ast))^2\frac{v_V^2}{a_V^2}\notag\\
        &\stackrel{(\ref{dpvh})}{=} \underbrace{\left(1 - (p_V^\ast h(p_V^\ast,p_L^\ast))^2\frac{v_V^2}{a_V^2}\right)}_{\geq 0}\left(\tau p_V^{\ast^2}h(p_V^\ast,p_L^\ast)[\![v]\!]^2 - 1\right)
        + \underbrace{2p_V^\ast(h(p_V^\ast,p_L^\ast))^2[\![v]\!]}_{< 0}.
    \end{align*}
    For the second term we obtain (as before in the proof of Lemma \ref{f_pL_pos})
    \begin{align*}
        &0 > \tau p_V^\ast(p_V^\ast h(p_V^\ast,p_L^\ast))[\![v]\!]^2 - 1\quad\stackrel{p_V^\ast h(p_V^\ast,p_L^\ast) > 0}{\Leftrightarrow}\quad
        \tau < \frac{1}{p_V^\ast(p_V^\ast h(p_V^\ast,p_L^\ast))[\![v]\!]^2}.
    \end{align*}
    and again we have
    \begin{align*}
        &\frac{1}{p_V^\ast(p_V^\ast h(p_V^\ast,p_L^\ast))[\![v]\!]^2} \geq \frac{v_V}{p_V^\ast a_V[\![v]\!]^2} = \frac{1}{p_V^\ast v_V a_V \left(\dfrac{v_L}{v_V} - 1\right)^2}
        \stackrel{(\ref{assumptions1})_1}{\geq} \frac{\gamma_V}{(1 - \alpha)^2a_V^3} \stackrel{(\ref{assumptions1})_3}{>} \tau.
    \end{align*}
    This proves the Lemma for $0 < p_V^\ast h(p_V^\ast,p_L^\ast) \leq a_V/v_V$.\\
    \newline
    \textbf{Third Case:} We discuss $a_V/v_V < p_V^\ast h(p_V^\ast,p_L^\ast) < a_L/v_L$ and rewrite $\partial_{p_V}f(p_V^\ast,p_L^\ast)$ to obtain with an analogue argument as used before
    \begin{align*}
        &\partial_{p_V}f(p_V^\ast,p_L^\ast) =\dots\\
        &= -1 + 2(p_V^\ast h(p_V^\ast,p_L^\ast))(h(p_V^\ast,p_L^\ast) + p_V^\ast\partial_{p_V}h(p_V^\ast,p_L^\ast))[\![v]\!] + (p_V^\ast h(p_V^\ast,p_L^\ast))^2\frac{v_V^2}{a_V^2}\\
        &= -\left(1 - (p_V^\ast h(p_V^\ast,p_L^\ast))^2\frac{v_V^2}{a_V^2}\right) + 2(p_V^\ast h(p_V^\ast,p_L^\ast))(p_V^\ast\partial_{p_V^\ast}h(p_V^\ast,p_L^\ast))[\![v]\!]
        + 2p_V^\ast h(p_V^\ast,p_l^\ast)^2[\![v]\!]\\
        &\stackrel{(\ref{dpvh})}{=} -\frac{2}{\tau [\![v]\!]}\partial_{p_V}h(p_V^\ast,p_L^\ast) + p_V^\ast h(p_V^\ast,p_L^\ast))(p_V^\ast\partial_{p_V^\ast}h(p_V^\ast,p_L^\ast))[\![v]\!]
        + 2p_V^\ast h(p_V^\ast,p_l^\ast)^2[\![v]\!]\\
        &= \underbrace{-\frac{2}{\tau [\![v]\!]}\partial_{p_V}h(p_V^\ast,p_L^\ast)}_{> 0}
           \underbrace{\left(1 - \tau p_V^\ast [\![v]\!]^2 p_V^\ast h(p_V^\ast,p_L^\ast)\right)}_{\stackrel{(\ref{assumptions1})}{< 0}}
         + \underbrace{2p_V^\ast h(p_V^\ast,p_l^\ast)^2[\![v]\!]}_{< 0} < 0.
    \end{align*}
    This ends the proof.
\end{prf}
\begin{cor}[Monotonicity of the Implicit Function]\label{impl_fun_mon}
    Let $(p_V^\ast,p_L^\ast)$ be an admissible solution $f(p_V^\ast,p_L^\ast) = 0$.
    Then there exists a function $\varphi$ with $p_L^\ast = \varphi(p_V^\ast)$ which is strictly monotonically increasing, i.e. $\varphi'(p_V^\ast) > 0$.
\end{cor}
\begin{prf}
    This follows using the implicit function theorem together with Lemma \ref{f_pL_pos} and \ref{f_pV_neg}.
\end{prf}
\begin{cor}\label{cond_evap_p_ineq}
    During a condensation process both pressures are larger than the saturation pressure
    \begin{align*}
        p_0 < p_V < p_L
    \end{align*}
    whereas during evaporation both pressures are smaller than the saturation pressure
    \begin{align*}
        p_V < p_L < p_0.
    \end{align*}
\end{cor}
\begin{prf}
    This follows from Corollary \ref{impl_fun_mon} and $p_L(p_0) = p_0$.
\end{prf}
\begin{lem}[Uniqueness]\label{unique_sol}
    Let $(p_V^\ast,p_L^\ast)$ be an admissible solution of $f(p_V^\ast,p_L^\ast) = 0$. Then this root is unique in the sense that for a given $p_V^\ast$ the solution $p_L^\ast$ is unique.
\end{lem}
\begin{prf}
    First we assume that there exists a $p_L^{\ast\ast} > p_L^\ast$ such that $f(p_V^\ast,p_L^{\ast\ast}) = 0$. From Lemma \ref{f_pL_pos} we know that $\partial_{p_L}f(p_V^\ast,p_L^\ast) > 0$.
    Hence we have (monotonicity argument) $\partial_{p_L}f(p_V^\ast,p_L^{\ast\ast}) \leq 0$. Therefore we have
    \begin{align}
        \underbrace{p_V^\ast h(p_V^\ast,p_L^{\ast\ast}) < -\frac{a_V}{v_V}}_{\mathbf{I}}
        \quad\dot{\vee}\quad
        \underbrace{p_V^\ast h(p_V^\ast,p_L^{\ast\ast}) > \frac{a_L}{v_L}}_{\mathbf{II}}\label{contradic_ineq}
    \end{align}
    otherwise we would meet the requirements of Lemma \ref{f_pL_pos}. Since $p_V^\ast h(p_V^\ast,p_L^\ast) \leq a_L/v_L$ and Lemma \ref{h_pL_decrease} we can exclude $\mathbf{II}$.
    Assuming $\mathbf{I}$ is true we have that the root $(p_V^\ast,p_L^{\ast\ast})$ itself is not admissible and every possible further root with $p_L > p_L^{\ast\ast}$ would also fulfill relation
    $\mathbf{I}$ due to Lemma \ref{h_pL_decrease} and thus is not admissible.\\
    Now we assume that there exists a $p_L^{\ast\ast} < p_L^\ast$ such that $f(p_V^\ast,p_L^{\ast\ast}) = 0$. As in the first case we have the two possibilities (\ref{contradic_ineq}).
    The arguments are now quite analogue to the first case. We can exclude $\mathbf{I}$ since
    \begin{align*}
        -\frac{a_V}{v_V} \leq p_V^\ast h(p_V^\ast,p_L^\ast) < p_V^\ast h(p_V^\ast,p_L^{\ast\ast}).
    \end{align*}
    Therefore relation $\mathbf{II}$ must hold and $p_L^{\ast\ast}$ is no admissible root. Due to Lemma \ref{h_pL_decrease} every further solution $p_L < p_L^{\ast\ast}$ also fulfills $\mathbf{II}$.
    This proves uniqueness.
\end{prf}
\begin{lem}[Global Existence]\label{glob_exist}
    For every $p_V^\ast \in [0,\sigma_{max}p_0]$ exists a $p_L^\ast \in [p_V^\ast,\infty)$ such that $(p_V^\ast,p_L^\ast)$ is an admissible root of $f(p_V^\ast,p_L^\ast) = 0$.
\end{lem}
\begin{prf}
    We already have local existence in a neighborhood of $(p_0,p_0)$ due to the implicit function theorem.
    In the following we discriminate the cases depending on whether $p_V$ is smaller or larger than the saturation pressure $p_0$.\\
    \newline
    \textbf{First Case $\mathbf{(0 \leq p_V < p_0)}$:}
    Assume that there exists a $p_V < p_0$ such that there exists \textit{no} $p_L$ with $f(p_V,p_L) = 0$.
    Using the above results we know that there exists an admissible root $(p_V^\ast,p_L^\ast)$ in the neighborhood of $(p_0,p_0)$
    and due to monotonicity/continuity a further root $p_V < \bar{p}_V < p_V^\ast$ and $\bar{p}_L$ such that
    \begin{align*}
        f(\bar{p}_V,\bar{p}_L) = 0\,\,\land\,\,\partial_{p_L}f(\bar{p}_V,\bar{p}_L) = 0.
    \end{align*}
    Hence this root is \textit{not} admissible due to Lemma \ref{f_pL_pos}. On the other hand we have, due to the behavior of the function $h(p_V,p_L)$ in $(p_0,p_0)$ and the fact that
    \begin{align*}
        f(p_V,p_L) = 0\,\,\land\,\,h(p_V,p_L) = 0\quad\Leftrightarrow\quad[\![p]\!] = 0,
    \end{align*}
    that $h(\bar{p}_V,\bar{p}_L) > 0$ for $\bar{p}_V < p_0$. Together with Corollary \ref{z_simple_upbound} this gives
    \begin{align*}
        0 < \bar{p}_V h(\bar{p}_V,\bar{p}_L) \leq \frac{a_L}{v_L}.
    \end{align*}
    This contradicts the above statement that the root $\bar{p}_V$ is not admissible. Therefore the nonexistence assumption is wrong and we have global existence for $0 \leq p_V < p_0$.\\
    \newline
    \textbf{Second Case $\mathbf{(p_0 < p_V \leq \sigma_{max}p_0)}$:} The idea is again to show, that there exists no $p_0 < p_V^\ast \leq \sigma_{max}p_0$ such that
    \begin{align}
        f(p_V^\ast,p_L^\ast) = 0\,\,\land\,\,\partial_{p_L}f(p_V^\ast,p_L^\ast) = 0.\label{contradict_cond}
    \end{align}
    Let us assume we have $(p_V^\ast, p_L^\ast)$ such that the above relation holds. From that we can conclude
    \begin{align*}
        &\partial_{p_L}f(p_V^\ast,p_L^\ast) = 0\quad\Leftrightarrow\\
        &\left(p_V^\ast h(p_V^\ast,p_L^\ast)\right)^2 = \left(1 + 2(p_V^\ast h(p_V^\ast,p_L^\ast))(p_V^\ast\partial_{p_L}h(p_V^\ast,p_L^\ast))[\![v]\!]\right)\frac{a_L^2}{v_L^2}.
    \end{align*}
    Inserting this expression in $0 = f(p_V^\ast,p_L^\ast)$ gives
    \begin{align*}
        0 = f(p_V^\ast,p_L^\ast) &= [\![p]\!] + \left(p_V^\ast h(p_V^\ast,p_L^\ast)\right)^2[\![v]\!]\\
        &= [\![p]\!] + \left(1 + 2(p_V^\ast h(p_V^\ast,p_L^\ast))(p_V^\ast\partial_{p_L}h(p_V^\ast,p_L^\ast))[\![v]\!]\right)\frac{a_L^2}{v_L^2}[\![v]\!]\\
        &= [\![p]\!] + [\![v]\!]\frac{a_L^2}{v_L^2} + 2(p_V^\ast h(p_V^\ast,p_L^\ast))(p_V^\ast\partial_{p_L}h(p_V^\ast,p_L^\ast))[\![v]\!]^2\\
        &= \frac{2}{\tau}\frac{a_L^2}{v_L^2}\partial_{p_L}h(p_V^\ast,p_L^\ast) + 2(p_V^\ast h(p_V^\ast,p_L^\ast))(p_V^\ast\partial_{p_L}h(p_V^\ast,p_L^\ast))[\![v]\!]^2\\
        &= \frac{2}{\tau}\partial_{p_L}h(p_V^\ast,p_L^\ast)\frac{a_L^2}{v_L^2}\left(1 + \tau \left.p_V^\ast\right.^2h(p_V^\ast,p_L^\ast)[\![v]\!]^2\right)
    \end{align*}
    We define the function
    \begin{align*}
        H(p_V,p_L) = 1 + \tau p_V^2h(p_V,p_L)[\![v]\!]^2.
    \end{align*}
    Due to Lemma \ref{h_pL_decrease} we have $H(p_V^\ast,p_L^\ast) = 0$ and hence
    \begin{align}
        p_V^\ast h(p_V^\ast,p_L^\ast) = -\frac{1}{\tau p_V^\ast[\![v]\!]^2}.\label{exist_cond_1}
    \end{align}
    Further we can rewrite $\partial_{p_L} f(p_V,p_L)$ in terms of $H(p_V,p_L)$, i.e.
    \begin{align*}
        \partial_{p_L}f(p_V,p_L) = -\frac{v_L^2}{\left(\tau p_Va_L[\![v]\!]^2\right)^2}(H(p_V,p_L) - 1)^2 + \left(1 + \frac{[\![p]\!]}{[\![v]\!]}\frac{v_L^2}{a_L^2}\right)(H(p_V,p_L) - 1) + 1.
    \end{align*}
    From this we immediately get
    \begin{align*}
        H(p_V,p_L) &= \partial_{p_L}f(p_V,p_L)
        \quad\Leftrightarrow\\
        0 &= (H(p_V,p_L) - 1)\left(-\frac{v_L^2}{\left(\tau p_Va_L[\![v]\!]^2\right)^2}(H(p_V,p_L) - 1) + \frac{[\![p]\!]}{[\![v]\!]}\frac{v_L^2}{a_L^2}\right).
    \end{align*}
    For the considered root $(p_V^\ast,p_L^\ast)$ we can exclude the first case since $H(p_V^\ast,p_L^\ast) = 1$ if and only if $p_V^\ast h(p_V^\ast,p_L^\ast) = 0$.
    Hence we further look at the second term which must vanish for $(p_V^\ast,p_L^\ast)$ and obtain
    \begin{align}
        H(p_V^\ast,p_L^\ast) = \left(\tau p_V^\ast[\![v]\!]^2\right)^2\frac{[\![p]\!]}{[\![v]\!]} + 1
        \quad\Leftrightarrow\quad
        p_V^\ast h(p_V^\ast,p_L^\ast) = \tau p_V^\ast[\![p]\!][\![v]\!].\label{exist_cond_2}
    \end{align}
    Summing up we can state that there are two conditions (\ref{exist_cond_1}) and (\ref{exist_cond_2}) which need to be true for $(p_V^\ast,p_L^\ast)$ when (\ref{contradict_cond}) holds.
    For equation (\ref{exist_cond_1}) we easily verify
    \begin{align}
        p_V^\ast h(p_V^\ast,p_L^\ast) = -\frac{1}{\tau p_V^\ast[\![v]\!]^2} \stackrel{(\ref{assumptions1})_1}{\leq} -\frac{1}{\tau p_V^\ast v_V^2(1 - \alpha)^2}
        \stackrel{(\ref{assumptions1})_3}{<} -\frac{a_V}{v_V}.\label{cond_1_upbound}
    \end{align}
    Now we investigate (\ref{exist_cond_2}) and prove that this implies $p_V^\ast h(p_V^\ast,p_L^\ast) > -a_V/v_V$. This would contradict (\ref{cond_1_upbound}) and hence finish the proof.\\
    First we introduce the following functions for fixed $p_V^\ast$
    \begin{align*}
        F(p_L) &:= p_V^\ast h(p_V^\ast,p_L), &F'(p_L) &= p_V^\ast\partial_{p_L}h(p_V^\ast,p_L),\\
        G(p_L) &:= \tau p_V^\ast[\![p]\!][\![v]\!], &G'(p_L) &= \tau p_V^\ast\left\{[\![v]\!] - [\![p]\!]\frac{v_L^2}{a_L^2}\right\}.
    \end{align*}
    We immediately verify for all $p_L \geq p_V^\ast$
    \begin{align*}
        G'(p_L) < F'(p_L) < 0.
    \end{align*}
    Furthermore we have
    \begin{align*}
        G(p_V^\ast) = 0 \stackrel{p_V^\ast > p_0}{>} F(p_V^\ast) = p_V^\ast h(p_V^\ast,p_V^\ast) = \tau p_V^\ast[\![g(p_V^\ast)]\!].
    \end{align*}
    Surely there is a $\bar{p}_L > p_V^\ast$ such that
    \begin{align}
        G(\bar{p}_L) = -\frac{a_V}{v_V}
        \quad\text{with}\quad
        \bar{p}_L = p_V^\ast - \frac{a_V}{\tau p_V^\ast v_V[\![v]\!]} \leq p_V^\ast + \frac{a_V}{\tau p_V^\ast v_V^2(1 - \alpha)}.\label{p_bar_ineq}
    \end{align}
    Now we investigate $F(\bar{p}_L)$ and obtain
    \begin{align*}
        F(\bar{p}_L) = p_V^\ast h(p_V^\ast,\bar{p}_L) &= \tau p_V^\ast\left\{[\![g]\!] - \frac{1}{2}[\![p]\!](v_L + v_V)\right\}
        \stackrel{(\ref{p_bar_ineq})}{=} \tau p_V^\ast([\![g]\!] - v_V[\![p]\!]) + \frac{1}{2}\frac{a_V}{v_V}.
    \end{align*}
    We have for $p_V^\ast = \sigma p_0$ with $\sigma \in [1,\sigma_{max}]$
    \begin{align*}
        &\tau p_V^\ast\left([\![g]\!] - v_V[\![p]\!]\right)\stackrel{g_L(\bar{p}_L) > g_L(p_V^\ast)}{>} \tau p_V^\ast\left([\![g(p_V^\ast)]\!] - v_V[\![p]\!]\right)\\
        &\stackrel{\substack{\hphantom{\frac{d}{dp}\frac{a(p)}{v(p)} > 0}\\{(\ref{p_bar_ineq})}}}{\geq} \tau p_V^\ast\left\{[\![g(p_V^\ast)]\!] - \frac{a_V}{\tau p_V^\ast v_V(1 - \alpha)}\right\}
        = \tau p_V^\ast [\![g(p_V^\ast)]\!] - \frac{a_V}{v_V(1 - \alpha)}\\
        &\stackrel{\substack{\hphantom{\frac{d}{dp}\frac{a(p)}{v(p)} > 0}\\ \text{Taylor}}}{\geq} \tau p_V^\ast [\![v(p_0)]\!](p_V^\ast - p_0) - \frac{a_V}{v_V(1 - \alpha)}
        > \frac{v_V(p_0)}{a_V(p_0)[\![v(p_0)]\!]}\sigma(\sigma - 1)p_0 - \frac{a_V}{v_V(1 - \alpha)}\\
        &\stackrel{\hphantom{\frac{d}{dp}\frac{a(p)}{v(p)} > 0}}{\geq} -\frac{a_V(p_0)}{v_V(p_0)}\frac{\sigma(\sigma - 1)}{1 - \alpha} - \frac{a_V}{v_V(1 - \alpha)}\\
        &\stackrel{\frac{\dd}{\dd p}\frac{a(p)}{v(p)} > 0}{>} -\frac{a_V(p_V^\ast)}{v_V(p_V^\ast)}\frac{\sigma(\sigma - 1) - 1}{1 - \alpha}\\
        &\stackrel{\substack{\hphantom{\frac{d}{dp}\frac{a(p)}{v(p)} > 0}\\ \sigma \leq \sigma_{max}}}{\geq} -\frac{3}{2}\frac{a_V}{v_V}.
    \end{align*}
    This gives us $F(\bar{p}_L) > G(\bar{p}_L)$ and so there exists a $p_L^{\ast\ast} \in (p_V^\ast, \bar{p}_L)$ such that
    \begin{align}
        F(p_L^{\ast\ast}) = G(p_L^{\ast\ast}) > -\frac{a_V}{v_V}.
    \end{align}
    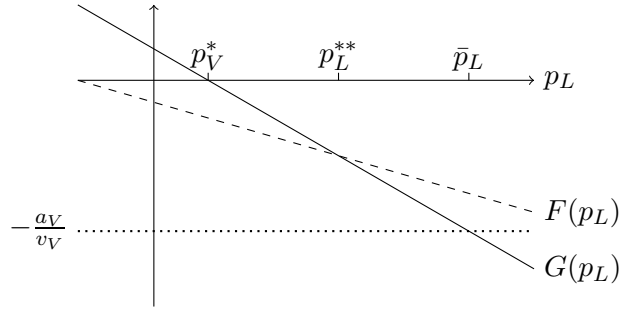
\begin{figure}[h!]
        \begin{center}
            \begin{tikzpicture}[domain=-4:4]
                %
                \draw[->] (-3,0) -- (3,0) node[right] {$p_L$};
                \draw (-9/7,0) node[above] {$p_V^\ast$};
                \draw (-9/7,0) -- (-9/7,0.1);
                \draw (3/7,0) node[above] {$p_L^{\ast\ast}$};
                \draw (3/7,0) -- (3/7,0.1);
                \draw (15/7,0) node[above] {$\bar{p}_L$};
                \draw (15/7,0) -- (15/7,0.1);
                \draw[->] (-2,-3) -- (-2,1);
                \draw[color=black,dotted,thick] (-3,-2) -- (3,-2);
                \draw (-3,-2) node[left] {$-\frac{a_V}{v_V}$};
                \draw[color=black,dashed] (-3,0) -- (3,-1.75);
                \draw (3,-1.75) node[right] {$F(p_L)$};
                \draw[color=black] (-3,1) -- (3,-2.5);
                \draw (3,-2.5) node[right] {$G(p_L)$};
            \end{tikzpicture}
        \end{center}
        \caption{Idea for the contradiction argument}
        \label{sketch_contradict_arg}
    \end{figure}
    Thus condition (\ref{exist_cond_2}) contradicts the first condition (\ref{exist_cond_1}). Hence there exists \textit{no} $(p_V^\ast,p_L^\ast)$ such that relation (\ref{contradict_cond}) holds.
    This implies global existence for all $p_0 < p_V \leq \sigma_{max}p_0$ and finishes the proof.
\end{prf}

    \newpage
    \section{Solution of the Two Phase Riemann Problem}\label{chap_2p_rp}
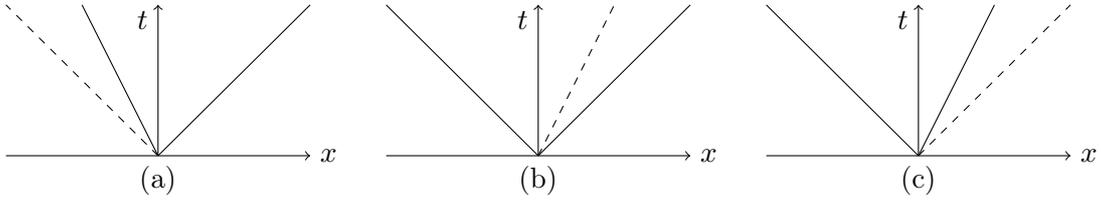
\begin{figure}[h!]
    \begin{center}
        \begin{tikzpicture}[domain=-7:7]
            %
            \draw[->] (-7,0) -- (-3,0) node[right] {$x$};
            \draw (-5,0) node[below] {(a)};
            \draw[->] (-5,0) -- (-5,2);
            \draw[color=black] (-5,1.8) node[left] {$t$};
            \draw[color=black,dashed] (-5,0) -- (-7,2);
            \draw[color=black] (-5,0) -- (-6,2);
            \draw[color=black] (-5,0) -- (-3,2);
            %
            \draw[->] (-2,0) -- (2,0) node[right] {$x$};
            \draw (0,0) node[below] {(b)};
            \draw[->] (0,0) -- (0,2);
            \draw[color=black] (0,1.8) node[left] {$t$};
            \draw[color=black] (0,0) -- (-2,2);
            \draw[color=black,dashed] (0,0) -- (1,2);
            \draw[color=black] (0,0) -- (2,2);
            %
            \draw[->] (3,0) -- (7,0) node[right] {$x$};
            \draw (5,0) node[below] {(c)};
            \draw[->] (5,0) -- (5,2);
            \draw[color=black] (5,1.8) node[left] {$t$};
            \draw[color=black] (5,0) -- (3,2);
            \draw[color=black] (5,0) -- (6,2);
            \draw[color=black,dashed] (5,0) -- (7,2);
        \end{tikzpicture}
    \end{center}
    \caption{Wave patterns. Solid line: classical waves. Dashed line: phase boundary}
    \label{3_diff_wave_patterns}
\end{figure}

\setcounter{para}{1}
In this section we want to solve the Riemann problem.
Therefore we follow the strategy of constructing \textit{wave curves} and obtain the solution as the intersection of the wave curves, as for example done in \cite{Menikoff1989,Toro2009}.
Due to the phase boundary we have an additional term, but we still want to show uniqueness of a solution to the Riemann problem.
Hence we need a further monotonicity argument which we will prove in the following.\\
To this end we additionally need bounds for the dimensionless speed of sound $\gamma_V$ and $\gamma_L$.
We distinguish two relevant cases, each with an appropriate condition needed to prove monotonicity.
This is necessary especially for EOS (or equivalently fluids) near the critical point, e.g.\ \textit{van der Waals} EOS.
Further these conditions show that the dimensionless quantities are not independent of each other. We consider the following relevant cases
\begin{align}
    (\textbf{I})&\begin{cases}
        \gamma_V &\leq 1\quad\text{and}\quad 1 \leq \gamma_L\quad\text{with},\\
        \dfrac{1}{\gamma_L} &\geq 1 + \dfrac{\varepsilon(\gamma_V)}{\alpha}.
        %
    \end{cases}\notag\\
    \label{assumptions2}\\
    (\textbf{II})&\begin{cases}
        \gamma_V &< 1\quad\text{and}\quad \gamma_L < 1\quad\text{with},\\
        \alpha &\leq \dfrac{1 - \frac{1}{\gamma_V}}{1 - \frac{1}{\gamma_L}}\quad\text{and}\quad\varepsilon(\gamma_V) \leq 0.
    \end{cases}\notag
\end{align}
The quantity $\varepsilon(\gamma_V)$ is defined as follows, using all quantities as introduced before,
\begin{align}
    \varepsilon(\gamma_V) := \frac{1}{\gamma_V} - 1 - \frac{\tau a_V^3}{\gamma_V^2}(1 - \alpha)^2\left(1 - (\alpha\beta)^2\right).\label{assum_const}
\end{align}
So far we proved in Section \ref{chap_pb_sol} that there exists a unique solution of the jump conditions at the interface.
Furthermore we can express the pressure in the liquid phase as a strictly monotone increasing function of the vapor pressure
\begin{align*}
    p_L^\ast = \varphi(p_V^\ast)\quad\text{with}\quad\varphi'(p_V^\ast) > 0.
\end{align*}
\begin{lem}\label{lem_mon_zvol}
    Given the requirements (\ref{assumptions1}) 
    and (\ref{assumptions2}). For an admissible solution $f(p_V^\ast,p_L^\ast) = 0$ the following monotonicity holds
    \begin{align*}
        \frac{\textup{d}}{\textup{d}p_V}\left(p_V^\ast h(p_V^\ast,p_L^\ast)[\![v]\!]\right) \geq 0.
    \end{align*}
\end{lem}
\begin{prf}
    We have
    \begin{align*}
        \frac{\textup{d}}{\textup{d}p_V}\left(p_V^\ast h(p_V^\ast,p_L^\ast)[\![v]\!]\right) &= \partial_{p_V}\left(p_V^\ast h(p_V^\ast,p_L^\ast)[\![v]\!]\right)
        + \partial_{p_L}\left(p_V^\ast h(p_V^\ast,p_L^\ast)[\![v]\!]\right)\varphi'(p_V^\ast)\\
        &= (h(p_V^\ast,p_L^\ast) + p_V^\ast\partial_{p_V}h(p_V^\ast,p_L^\ast))[\![v]\!] + p_V^\ast h(p_V^\ast,p_L^\ast)\frac{v_V^2}{a_V^2}\\
        &+ \left\{p_V^\ast\partial_{p_L}h(p_V^\ast,p_L^\ast)[\![v]\!] - p_V^\ast h(p_V^\ast,p_L^\ast)\frac{v_L^2}{a_L^2}\right\}\varphi'(p_V^\ast).\\
    \end{align*}
    For $(p_0,p_0)$ the statement is obvious and hence we assume $p_V^\ast h(p_V^\ast,p_L^\ast) \neq 0$ from now on. Now we can write
    \begin{align*}
        \frac{\textup{d}}{\textup{d}p_V}\left(p_V^\ast h(p_V^\ast,p_L^\ast)[\![v]\!]\right)
        &= \frac{1}{2p_V^\ast h(p_V^\ast,p_L^\ast)}(\partial_{p_V}f(p_V^\ast,p_L^\ast) + 1) + \frac{1}{2}p_V^\ast h(p_V^\ast,p_L^\ast)\frac{v_V^2}{a_V^2}\\
        &+ \left\{\frac{1}{2p_V^\ast h(p_V^\ast,p_L^\ast)}(\partial_{p_L}f(p_V^\ast,p_L^\ast) - 1) - \frac{1}{2}p_V^\ast h(p_V^\ast,p_L^\ast)\frac{v_L^2}{a_L^2}\right\}\varphi'(p_V^\ast).\\
    \end{align*}
    We multiply with $\partial_{p_L}f(p_V^\ast,p_L^\ast) > 0$ and use
    \begin{align*}
        \varphi'(p_V^\ast) = -\frac{\partial_{p_V}f(p_V^\ast,p_L^\ast)}{\partial_{p_L}f(p_V^\ast,p_L^\ast)}.
    \end{align*}
    Thus we obtain
    \begin{align*}
        &\partial_{p_L}f(p_V^\ast,p_L^\ast)\frac{\textup{d}}{\textup{d}p_V}\left(p_V^\ast h(p_V^\ast,p_L^\ast)[\![v]\!]\right) = \dots\\
        &=\left\{\frac{1}{2p_V^\ast h(p_V^\ast,p_L^\ast)}(\partial_{p_V}f(p_V^\ast,p_L^\ast) + 1)
        + \frac{1}{2}p_V^\ast h(p_V^\ast,p_L^\ast)\frac{v_V^2}{a_V^2}\right\}\partial_{p_L}f(p_V^\ast,p_L^\ast)\\
        &- \left\{\frac{1}{2p_V^\ast h(p_V^\ast,p_L^\ast)}(\partial_{p_L}f(p_V^\ast,p_L^\ast) - 1)
        - \frac{1}{2}p_V^\ast h(p_V^\ast,p_L^\ast)\frac{v_L^2}{a_L^2}\right\}\partial_{p_V}f(p_V^\ast,p_L^\ast)\\
        &= \frac{1}{2p_V^\ast h(p_V^\ast,p_L^\ast)}(\partial_{p_V}f(p_V^\ast,p_L^\ast) + \partial_{p_L}f(p_V^\ast,p_L^\ast))\\
        &+ \frac{1}{2}p_V^\ast h(p_V^\ast,p_L^\ast)\left(\partial_{p_V}f(p_V^\ast,p_L^\ast)\frac{v_L^2}{a_L^2} + \partial_{p_L}f(p_V^\ast,p_L^\ast)\frac{v_V^2}{a_V^2}\right)\\
        &= (h(p_V^\ast,p_L^\ast) + p_V\partial_{p_V}h(p_V^\ast,p_L^\ast))\left(1 + (p_V^\ast h(p_V^\ast,p_L^\ast))^2\frac{v_L^2}{a_L^2}\right)[\![v]\!]\\
        &+ p_V\partial_{p_L}h(p_V^\ast,p_L^\ast)\left(1 + (p_V^\ast h(p_V^\ast,p_L^\ast))^2\frac{v_V^2}{a_V^2}\right)[\![v]\!]
        + p_V^\ast h(p_V^\ast,p_L^\ast)\left(\frac{v_V^2}{a_V^2} - \frac{v_L^2}{a_L^2}\right)\\
        &= h(p_V^\ast,p_L^\ast)[\![v]\!]\left(1 + (p_V^\ast h(p_V^\ast,p_L^\ast))^2\frac{v_L^2}{a_L^2}\right)
        - p_V^\ast h(p_V^\ast,p_L^\ast)[\![\frac{v^2}{a^2}]\!]\\
        &+ \tau p_V[\![v]\!]^2\left(1 - (p_V^\ast h(p_V^\ast,p_L^\ast))^4\frac{v_V^2}{a_V^2}\frac{v_L^2}{a_L^2}\right)\\
        &= h(p_V^\ast,p_L^\ast)\left([\![v]\!] + [\![v]\!](p_V^\ast h(p_V^\ast,p_L^\ast))^2\frac{v_L^2}{a_L^2} - p_V^\ast[\![\frac{v^2}{a^2}]\!]\right)\\
        &+ \tau p_V[\![v]\!]^2\left(1 - (p_V^\ast h(p_V^\ast,p_L^\ast))^4\frac{v_V^2}{a_V^2}\frac{v_L^2}{a_L^2}\right)\\
        &= h(p_V^\ast,p_L^\ast)\left(v_L\left(1 - \frac{p_L^\ast v_L}{a_L^2}\right) - v_V\left(1 - \frac{p_V^\ast v_V}{a_V^2}\right)\right)\\
        &+ \tau p_V[\![v]\!]^2\left(1 - (p_V^\ast h(p_V^\ast,p_L^\ast))^4\frac{v_V^2}{a_V^2}\frac{v_L^2}{a_L^2}\right)
    \end{align*}
    Due to the bounds for the \textit{EOS} we can show
    \begin{align*}
        &\left(v_L\left(1 - \frac{p_L^\ast v_L}{a_L^2}\right) - v_V\left(1 - \frac{p_V^\ast v_V}{a_V^2}\right)\right) \geq 0.
        %
    \end{align*}
    and hence we can immediately verify the Lemma for
    \begin{align}
        0 < p_V^\ast h(p_V^\ast,p_L^\ast) \leq \sqrt{\dfrac{a_V}{v_V}\dfrac{a_L}{v_L}}.
    \end{align}
    Now we want to prove the result for $0 > p_V^\ast h(p_V^\ast,p_L^\ast) \geq -a_V/v_V$. We have
    \begin{align*}
        &h(p_V^\ast,p_L^\ast)\left(v_L\left(1 - \frac{p_L^\ast v_L}{a_L^2}\right) - v_V\left(1 - \frac{p_V^\ast v_V}{a_V^2}\right)\right)
        + \tau p_V^\ast[\![v]\!]^2\left(1 - (p_V^\ast h(p_V^\ast,p_L^\ast))^4\frac{v_V^2}{a_V^2}\frac{v_L^2}{a_L^2}\right)\\
        &\geq -\frac{a_V}{p_V^\ast v_V}\left(v_L\left(1 - \frac{1}{\gamma_L}\right) - v_V\left(1 - \frac{1}{\gamma_V}\right)\right)
        + \tau p_V^\ast[\![v]\!]^2\left(1 - \frac{a_V^2}{v_V^2}\frac{v_L^2}{a_L^2}\right)\\
        &\geq -\frac{\gamma_V}{a_V}\left(v_L\left(1 - \frac{1}{\gamma_L}\right) - v_V\left(1 - \frac{1}{\gamma_V}\right)\right) + \tau v_V^2p_V^\ast(1 - \alpha)^2\left(1 - (\alpha\beta)^2\right)\\
        &= -\frac{v_V}{a_V}\left(\frac{\gamma_V}{\gamma_L}\frac{v_L}{v_V}\left(\gamma_L - 1\right) + \left(1 - \gamma_V\right)\right)
        + \tau v_V^2p_V^\ast(1 - \alpha)^2\left(1 - (\alpha\beta)^2\right)\\
        &= -\frac{v_V}{a_V}\left(\frac{\gamma_V}{\gamma_L}\frac{v_L}{v_V}\left(\gamma_L - 1\right) + \left(1 - \gamma_V\right)
        - \frac{\tau a_V^3}{\gamma_V}(1 - \alpha)^2\left(1 - (\alpha\beta)^2\right)\right)\tag{+}\\
        &\stackrel{\gamma_L \geq 1}{\geq} -\frac{v_V}{a_V}\left(\alpha\frac{\gamma_V}{\gamma_L}\left(\gamma_L - 1\right) + \left(1 - \gamma_V\right)
        - \frac{\tau a_V^3}{\gamma_V}(1 - \alpha)^2\left(1 - (\alpha\beta)^2\right)\right)\\
        &\geq \alpha\gamma_V\frac{v_V}{a_V}\left(\frac{1}{\gamma_L} - \left(1 + \frac{\varepsilon(\gamma_V)}{\alpha}\right)\right)\\
        &\stackrel{(\ref{assumptions2}) (\textbf{I})}{\geq} 0.
    \end{align*}
    Starting from $(+)$ we obtain for the case $1 > \gamma$ for both phases
    \begin{align*}
        &-\frac{v_V}{a_V}\left(\frac{\gamma_V}{\gamma_L}\frac{v_L}{v_V}\left(\gamma_L - 1\right) + \left(1 - \gamma_V\right)
        - \frac{\tau a_V^3}{\gamma_V}(1 - \alpha)^2\left(1 - (\alpha\beta)^2\right)\right)\\
        &\geq -\frac{v_V}{a_V}\left(1 - \gamma_V - \frac{\tau a_V^3}{\gamma_V}(1 - \alpha)^2\left(1 - (\alpha\beta)^2\right)\right)\\
        &= -\gamma_V\frac{v_V}{a_V}\varepsilon(\gamma_V)\\
        &\stackrel{(\ref{assumptions2}) (\textbf{II})}{\geq} 0.
    \end{align*}
    It remains the case for
    \begin{align*}
        \sqrt{\frac{a_V}{v_V}\frac{a_L}{v_L}} < p_V^\ast h(p_V^\ast,p_L^\ast) < \frac{a_L}{v_L}.
    \end{align*}
    In the subsequent Lemma \ref{z_max} we will exclude this case and thus the proof of this Lemma is finished.
\end{prf}\\
\begin{rem}[Assumptions on $\gamma$]
    In Lemma \ref{lem_mon_zvol} we only consider cases where $\gamma_V \in (0,1]$ for the vapor phase. As mentioned before the lower bound ensures hyperbolicity and thermodynamic stability.
    The upper bound is due to the fact, that we only consider pressures and temperatures \textup{below} the critical point.\\
    To illustrate this we consider the \textup{isothermal compressibility} $\kappa_T$ which is defined as follows
    \begin{align*}
        \kappa_T = -v\left(\frac{\partial v}{\partial p}\right)_T.
    \end{align*}
    For \textup{real} gases $\kappa_T$ can be expressed in terms of the pressure and the \textup{compressibility} or \textup{gas deviation factor} $\mathcal{Z}$
    (not to confuse with the mass flux $z$ used in this work), i.e.
    \begin{align*}
        \kappa_T = \frac{1}{p} - \frac{1}{\mathcal{Z}}\left(\frac{\partial \mathcal{Z}}{\partial p}\right)_T.
    \end{align*}
    Below the critical point the second term is negative for most gases and hence
    \begin{align*}
        \kappa_T > \frac{1}{p}\quad\Leftrightarrow\quad \gamma_V = \frac{1}{p\kappa_T} < 1.
    \end{align*}
    This property is reflected by nonlinear EOS such as the \textup{van der Waals} or \textup{Dieterici} EOS.
    For an ideal gas the second term vanishes and we obtain $\gamma_V = 1$.
\end{rem}
\begin{lem}\label{z_max}
    Consider two phases such that the requirements (\ref{assumptions1}) are fulfilled.
    Then there exists a maximal mass flux $z_{max}$ such that for every admissible solution $f(p_V^\ast,p_L^\ast) = 0$ the following upper bound holds
    \begin{align*}
        z_{max} \leq \sqrt{\frac{a_V}{v_V}\frac{a_L}{v_L}}.
    \end{align*}
\end{lem}
\begin{prf}
    Since $z(p_V) = 0$ if and only if $p_V = p_0$ and further $z(p_0)' < 0$ we can focus on vapor pressures smaller than $p_0$. We assume that
    \begin{align*}
        z_{max} > \sqrt{\frac{a_V}{v_V}\frac{a_L}{v_L}}.
    \end{align*}
    Hence there exists a $\tilde{p} \in (0,p_0)$ such that
    \begin{align*}
        z(\tilde{p}) = \tilde{p}h(\tilde{p},\varphi(\tilde{p})) = \sqrt{\frac{a_V}{v_V}\frac{a_L}{v_L}}\,\,\text{and}\,\,z'(\tilde{p}) \leq 0.
    \end{align*}
    This gives
    \begin{align*}
        0 &\geq z'(\tilde{p}) = \frac{1}{\tilde{p}}\sqrt{\frac{a_V}{v_V}\frac{a_L}{v_L}}
        + \tilde{p}\left(\partial_{p_V}h(\tilde{p},\varphi(\tilde{p})) + \partial_{p_L}h(\tilde{p},\varphi(\tilde{p}))\varphi'(\tilde{p})\right)\\
        &= \frac{1}{\tilde{p}}\sqrt{\frac{a_V}{v_V}\frac{a_L}{v_L}}
        + \frac{\tau\tilde{p}}{2}[\![v]\!]\left(1 - \frac{a_L}{a_V}\frac{v_V}{v_L} + \left(1 - \frac{a_V}{a_L}\frac{v_L}{v_V}\right)\varphi'(\tilde{p})\right)\\
        &= \underbrace{\frac{1}{\tilde{p}}\sqrt{\frac{a_V}{v_V}\frac{a_L}{v_L}}}_{> 0}
        + \underbrace{\frac{\tau\tilde{p}}{2}[\![v]\!]\left(1 - \frac{a_V}{a_L}\frac{v_L}{v_V}\right)}_{< 0}\left(\varphi'(\tilde{p}) - \frac{a_L}{a_V}\frac{v_V}{v_L}\right)\\
        &\Rightarrow \xi := \varphi'(\tilde{p}) > \frac{a_L}{a_V}\frac{v_V}{v_L} \geq \frac{1}{\alpha\beta} \stackrel{(\ref{assumptions1})_2}{>} 1.
    \end{align*}
    Using the definition of $\varphi'(p_V)$ we obtain
    \begin{align}
        &-\partial_{p_V}f(\tilde{p},\varphi(\tilde{p})) = \partial_{p_L}f(\tilde{p},\varphi(\tilde{p}))\xi\notag\\
        &\Leftrightarrow\notag\\
        &1 - 2\sqrt{\frac{a_V}{v_V}\frac{a_L}{v_L}}\left(\frac{1}{\tilde{p}}\sqrt{\frac{a_V}{v_V}\frac{a_L}{v_L}}
        + \tilde{p}\partial_{p_V}h(\tilde{p},\varphi(\tilde{p}))\right)[\![v]\!] - \frac{a_L}{v_L}\frac{v_V}{a_V}\notag\\
        &= \xi\left(1 + 2\sqrt{\frac{a_V}{v_V}\frac{a_L}{v_L}}\tilde{p}\partial_{p_L}h(\tilde{p},\varphi(\tilde{p}))[\![v]\!] - \frac{a_V}{v_V}\frac{v_L}{a_L}\right)\notag\\
        &\Leftrightarrow\notag\\
        &1 - \xi + \frac{a_V}{v_V}\frac{a_L}{v_L}\left(\xi\frac{v_L^2}{a_L^2} - \frac{v_V^2}{a_V^2}\right) = \dots\notag\\
        &= 2\sqrt{\frac{a_V}{v_V}\frac{a_L}{v_L}}[\![v]\!]\left(\frac{1}{\tilde{p}}\sqrt{\frac{a_V}{v_V}\frac{a_L}{v_L}} + \tilde{p}(\partial_{p_V}h(\tilde{p},\varphi(\tilde{p}))
        + \xi\partial_{p_L}h(\tilde{p},\varphi(\tilde{p})))\right).\label{z_max_contradict_eq}
    \end{align}
    For the right hand side of (\ref{z_max_contradict_eq}) we easily see
    \begin{align*}
        2\sqrt{\frac{a_V}{v_V}\frac{a_L}{v_L}}[\![v]\!]\left(\frac{1}{\tilde{p}}\sqrt{\frac{a_V}{v_V}\frac{a_L}{v_L}} + \tilde{p}(\partial_{p_V}h + \xi\partial_{p_L}h)\right)
        = 2\sqrt{\frac{a_V}{v_V}\frac{a_L}{v_L}}[\![v]\!]z'(\tilde{p}) \geq 0.
    \end{align*}
    If we consider the left hand side of (\ref{z_max_contradict_eq}) as a function of $\xi$ we get
    \begin{align*}
        \frac{\textup{d}}{\textup{d}\xi}\left(1 - \xi + \frac{a_V}{v_V}\frac{a_L}{v_L}\left(\xi\frac{v_L^2}{a_L^2} - \frac{v_V^2}{a_V^2}\right)\right)
        = -1 + \frac{a_V}{a_L}\frac{v_L}{v_V} \leq -1 + \alpha\beta < 0.
    \end{align*}
    Thus the left hand side of (\ref{z_max_contradict_eq}) is strictly decreasing in $\xi$ and we have
    \begin{align*}
        1 - \xi + \frac{a_V}{v_V}\frac{a_L}{v_L}\left(\xi\frac{v_L^2}{a_L^2} - \frac{v_V^2}{a_V^2}\right)
        \stackrel{\xi = 1}{=} \frac{a_V}{v_V}\frac{a_L}{v_L}[\![\frac{v^2}{c^2}]\!] < 0.
    \end{align*}
    Since $\xi > 1$ the left hand side of (\ref{z_max_contradict_eq}) is negative and hence contradicts the positive right hand side. Therefore the assumption for $z_{max}$ is wrong.
\end{prf}
\begin{rem}
    Lemma \ref{z_max} is a direct improvement of Corollary \ref{z_simple_upbound} obtained during the proof of Theorem \ref{exis_uniq_thm}.
    There we stated that the upper bound $a_L/v_L$ for the mass flux $z$ is always fulfilled.
\end{rem}
Now we consider two phase flows, where we initially have the vapor phase on the left ($x < 0$) and the liquid phase on the right side ($x > 0$).
The different phases are described using the corresponding \textit{EOS}. The considered Riemann initial data is
\begin{align}
    \rho(x,0) = \begin{cases} &\rho_V,\,x < 0\\ &\rho_L,\,x > 0\end{cases}
    \quad\text{and}\quad
    u(x,0) = \begin{cases} &u_V,\,x < 0\\ &u_L,\,x > 0\end{cases}.\label{init_data_rp}
\end{align}
The solution consists of two classical waves and the phase boundary separating four constant states. Hence there are three possible wave patterns, see Figure \ref{3_diff_wave_patterns}.
\subsection[First Case: Two Phase Flow without Phase Transition]{$1^{st}$ Case: Two Phase Flow without Phase Transition}\label{2p_wo_pt}
At first we want to deal with the case where phase transition is excluded, i.e. $z = 0$. Let us consider a wave pattern of type (b), see Figure \ref{3_diff_wave_patterns}.
The four constant states are denoted as follows
\begin{align*}
    \mathbf{W}_V = \left(\begin{matrix} \rho_V\\ u_V \end{matrix}\right),\quad
    \mathbf{W}_V^\ast = \left(\begin{matrix} \rho_V^\ast\\ u_V^\ast \end{matrix}\right),\quad
    \mathbf{W}_L^\ast = \left(\begin{matrix} \rho_L^\ast\\ u_L^\ast \end{matrix}\right),\quad
    \mathbf{W}_L = \left(\begin{matrix} \rho_L\\ u_L \end{matrix}\right).
\end{align*}
As in Section \ref{sec_riemann_problem} we want to derive a single function such that the single root $p$ is the solution for the pressure $p_V^\ast$. This procedure again uses the constancy of
pressure and velocity across the phase boundary, $u_V^\ast = u_L^\ast$ and $p_V^\ast = p_L^\ast$, which is because of $z = 0$.
For the solution we use the results obtained in Section \ref{chap_euler_eq}.
\begin{thm}[Solution without Phase Transition]\label{sol_without_pt}
    Let $f(p,\mathbf{W}_V,\mathbf{W}_L)$ be given as
    \begin{align*}
        f(p,\mathbf{W}_V,\mathbf{W}_L) = f_V(p,\mathbf{W}_V) + f_L(p,\mathbf{W}_L) + \Delta u,\,\,\Delta u = u_L - u_V,
    \end{align*}
    with the functions $f_V$ and$f_L$ given by
    \begin{align*}
        f_K(p,\mathbf{W}_K) &=
        \begin{dcases}
            \sqrt{-[\![p]\!][\![v_K]\!]},\,\,p > p_K\,\,\text{(Shock)}\\
            \int_{p_K}^p \frac{v_K(\zeta)}{a_K(\zeta)}\,\textup{d}\zeta,\,\,p \leq p_K\,\,\text{(Rarefaction)}
        \end{dcases},\quad K \in \{V,L\}.
    \end{align*}
    If there is a root $f(p^\ast,\mathbf{W}_V,\mathbf{W}_L) = 0$ with $0 < p^\ast \leq \tilde{p}$, then this root is unique. Here $\tilde{p}$ is given as in Definition \ref{max_pV}.
    Further this is the unique solution for the pressure $p_V^\ast$ of the Riemann problem (\ref{mass_cons})-(\ref{mom_balance}), (\ref{init_data_rp}).
    The velocity $u^\ast := u_V^\ast = u_L^\ast$ is given by
    \begin{align*}
        u^\ast = \frac{1}{2}(u_L + u_V) + \frac{1}{2}(f_L(p^\ast,\mathbf{W}_V) - f_V(p^\ast,\mathbf{W}_L)).
    \end{align*}
\end{thm}
\begin{prf}
    The function $f$ is strictly monotone increasing in $p$ due to the inequalities (\ref{mono_conc_shock}), (\ref{mono_conc_raref}) and Lemma \ref{lem_mon_zvol}.
    Furthermore we have $f(p,\mathbf{W}_V,\mathbf{W}_L) \to -\infty$ for $p \to 0$. Hence $f$ has at most one unique root, which is by construction the solution for the pressure $p_V^\ast$.
    The statement for the velocity $u^\ast$ follows immediately from the results in Section \ref{chap_euler_eq}.
\end{prf}
\newline
Note that one has to choose the corresponding \textit{EOS} to calculate the pressure depending quantities according to the index $K \in \{L,V\}$.
\begin{thm}[Sufficient Condition for Solvability]\label{suff_cond_solv_0}
    Consider the Riemann problem (\ref{mass_cons})-(\ref{mom_balance}), (\ref{init_data_rp}). We have two cases.
    \begin{itemize}
        \item[(i)] For $p_L < p_V = \tilde{p}$ the considered Riemann problem is solvable if and only if
        \begin{align*}
            &f(\tilde{p},\mathbf{W}_V,\mathbf{W}_L) = \dots\\
            &= \sqrt{-(\tilde{p} - p_V)(v_V(\tilde{p}) - v_V(p_V))} + \sqrt{-(p_L - \tilde{p})(v_L(p_L) - v_L(\tilde{p}))} + \Delta u \geq 0.
        \end{align*}
        \item[(ii)] For $p_L \geq \tilde{p}$ the considered Riemann problem is solvable if and only if
        \begin{align*}
            f(\tilde{p},\mathbf{W}_V,\mathbf{W}_L) = \sqrt{-(\tilde{p} - p_V)(v_V(\tilde{p}) - v_V(p_V))} + \int_{p_L}^{\tilde{p}} \frac{v_L(\zeta)}{a_L(\zeta)}\,\textup{d}\zeta + \Delta u \geq 0.
        \end{align*}
    \end{itemize}
\end{thm}
\begin{prf}
    As seen before in the proof of Theorem \ref{sol_without_pt}, $f$ is strictly monotone increasing in $p$ with $f(p,\mathbf{W}_V,\mathbf{W}_L) \to -\infty$ for $p \to 0$.
    Accordingly $f$ has a unique root if and only if $f(p,\mathbf{W}_V,\mathbf{W}_L) \geq 0$ for $p \to \tilde{p}$.
\end{prf}
\newline
So far we discussed the case that the solution is of type (b). The following result deals with the cases (a) and (c).
\begin{lem}\label{exclude_wave_patterns}
    There exists no solution of wave pattern types (a) and (c). This includes the coincidence of a classical wave and the phase boundary.
\end{lem}
\begin{prf}
    Let us first discuss case (c). For the notation see Figure \ref{wave_pattern_c}.
    \begin{figure}[h!]
    \begin{center}
        \begin{tikzpicture}[domain=-4:4]
            %
            \draw[->] (-4,0) -- (4,0) node[right] {$x$};
            \draw[->] (0,0) -- (0,4);
            \draw[color=black] (0,3.8) node[left] {$t$};
            \draw (-2.5,3) node[left] {$\mathbf{W}_V$};
            \draw[color=black] (0,0) -- (-3,4);
            \draw (-0.5,3) node[left] {$\mathbf{W}_V^{\ast\ast}$};
            \draw[color=black] (0,0) -- (1,4);
            \draw (1.85,3) node[left] {$\mathbf{W}_V^\ast$};
            \draw (3.75,3) node[left] {$\mathbf{W}_L$};
            \draw[color=black,dashed] (0,0) -- (3,4);
        \end{tikzpicture}
    \end{center}
    \caption{Wave pattern of type (c) with notation}
    \label{wave_pattern_c}
    \end{figure}
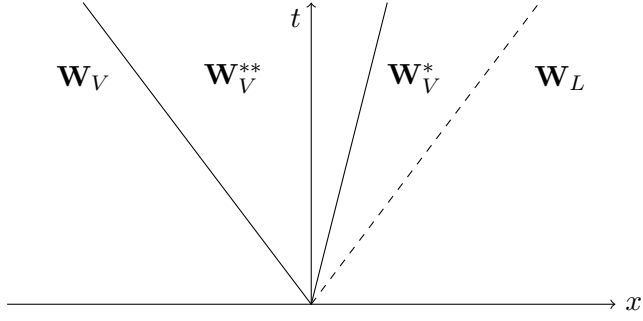
    Since $z = 0$ we have $w = u_L = u_V^\ast$ for the velocity of the interface. Further we assume, that the right classical wave is a shock moving with speed $S$.
    It is obvious that $w \geq S$ must hold.
    For the case of a right shock we have $p_V^{\ast\ast} > p_V^\ast$ and hence we obtain from the entropy inequality $Q > 0$, see Subsection \ref{sec_entropy_ineq}.
    Now we make use of the continuity of the mass flux across a shock wave and obtain
    \begin{align*}
        Q = -\rho_V^\ast(u_V^\ast - S)\quad\Leftrightarrow\quad u_V^\ast - S = -\frac{Q}{\rho_V^\ast}\quad\stackrel{u_V^\ast = w}{\Leftrightarrow}\quad S = w + \frac{Q}{\rho_V^\ast} > w.
    \end{align*}
    This contradicts the condition $w \geq S$.
    If, on the other hand, the right classical wave is a rarefaction wave we have for the head speed $u_V^\ast + a_V(p_V^\ast)$, see Subsection \ref{sec_rarefaction_wave}.
    Again this contradicts $w = u_L = u_V^\ast \geq u_V^\ast + a_V(p_V^\ast)$. In case that the phase boundary lies inside the rarefaction wave, we obtain similar contradictions in the wave speeds.
    For wave pattern type (a) the arguments are analogue.
\end{prf}
\subsection[Second Case: Two Phase Flow with Phase Transition]{$2^{nd}$ Case: Two Phase Flow with Phase Transition}\label{2p_w_pt}
Now we want to take phase transition into account, i.e. $z \neq 0$. As before we first want to discuss the wave pattern of type (b), see Figure \ref{3_diff_wave_patterns}.
In order to determine the solution we again construct a function analogue to Subsection \ref{2p_wo_pt}. For the left and right classical waves we use
\begin{align}
    u_V^\ast = u_V - f_V(p_V^\ast,\mathbf{W}_V)\quad\text{and}\quad u_L^\ast = u_L + f_L(p_L^\ast,\mathbf{W}_L).\label{star_state_velocities}
\end{align}
Across the phase boundary we make use of the jump conditions and obtain as in Subsection \ref{sec_entropy_ineq}
\begin{align}
    [\![u]\!] = u_L^\ast - u_V^\ast = -z[\![v]\!] = -z(v_L(p_L^\ast) - v_V(p_V^\ast)).\label{jump_velocities_z}
\end{align}
Finally we use the results obtained in Section \ref{chap_pb_sol}, especially Theorem \ref{exis_uniq_thm},
to express the liquid pressure at the interface as a function of the interface vapor pressure $p_L = \varphi(p_V)$. Combining these considerations we end up with the following theorem.
\begin{thm}[Solution with Phase Transition]\label{sol_with_pt}
    Let $f_z(p,\mathbf{W}_V,\mathbf{W}_L)$ be given as
    \begin{align*}
        f_z(p,\mathbf{W}_V,\mathbf{W}_L) = f_V(p,\mathbf{W}_V) + f_L(\varphi(p),\mathbf{W}_L) + z[\![v]\!] + \Delta u,\,\,\Delta u = u_L - u_V,
    \end{align*}
    with the functions $f_V$ and$f_L$ given by
    \begin{align*}
        f_V(p,\mathbf{W}_V) &=
        \begin{dcases}
            \sqrt{-[\![p]\!][\![v]\!]},\,\,p > p_V\,\,\text{(Shock)}\\
            \int_{p_V}^p \frac{v_V(\zeta)}{a_V(\zeta)}\,\textup{d}\zeta,\,\,p \leq p_V\,\,\text{(Rarefaction)}
        \end{dcases},\\
        f_L(\varphi(p),\mathbf{W}_L) &=
        \begin{dcases}
            \sqrt{-[\![p]\!][\![v]\!]},\,\,\varphi(p) > p_L\,\,\text{(Shock)}\\
            \int_{p_L}^{\varphi(p)} \frac{v_L(\zeta)}{a_L(\zeta)}\,\textup{d}\zeta,\,\,\varphi(p) \leq p_L\,\,\text{(Rarefaction)}
        \end{dcases}.
    \end{align*}
    The function $\varphi(p)$ is implicitly defined by (\ref{pb_mom_balance}) and the mass flux is given by (\ref{kin_rel2}).
    If there is a root $f_z(p^\ast,\mathbf{W}_V,\mathbf{W}_L) = 0$ with $0 < p^\ast \leq \tilde{p}$, this root is unique. If further
    \begin{align}
        p^\ast > p_V\quad\text{we must have}\quad z > -\frac{a_V(\bar{p})}{v_V(\bar{p})}\quad\text{for}\quad\bar{p} \in (p_V, p^\ast).\label{vap_shock_z_bound}
    \end{align}
    Then $p^\ast$ is the unique solution for the pressure $p_V^\ast$ of a (b)-type solution of the Riemann problem (\ref{mass_cons})-(\ref{mom_balance}), (\ref{init_data_rp}) with phase transition.
    If there is no root or condition (\ref{vap_shock_z_bound}) is not satisfied, then there is no solution to the mentioned Riemann problem.
\end{thm}
\begin{prf}
    Due to (\ref{mono_conc_shock}), (\ref{mono_conc_raref}), Corollary \ref{impl_fun_mon} and Lemma \ref{lem_mon_zvol} we get that the function $f_z$ is strictly monotone increasing in $p$.
    Furthermore we have $f(p,\mathbf{W}_V,\mathbf{W}_L) \to -\infty$ for $p \to 0$. Hence $f$ has at most one unique root, which is by construction the solution for the pressure $p_V^\ast$.
    Theorem \ref{exis_uniq_thm} then uniquely defines the liquid pressure $p_L^\ast = \varphi(p_V^\ast)$ and the mass flux $z$ at the interface. The remaining quantities can be calculated using the
    \textit{EOS} and (\ref{star_state_velocities}).
\end{prf}
\begin{rem}
    Condition (\ref{vap_shock_z_bound}) is needed in the case of a shock wave in the vapor phase to guarantee that $w > S$.
    Where $w$ denotes the velocity of the interface and $S$ of the shock respectively. This can be obtained as follows
    \begin{align*}
        &u_V^\ast - S = -v_V(p_V^\ast)Q_S\quad\text{and}\quad u_V^\ast - w = -v_V(p_V^\ast)z\quad\Leftrightarrow\quad w - S = v_V(p_V^\ast)(z - Q_S)\\
        &\Rightarrow\quad w > S\quad\Leftrightarrow\quad z > Q_S = - \frac{a_V(\bar{p})}{v_V(\bar{p})}.
    \end{align*}
    For the last equality we used the Lax condition for $S$ together with the monotonicity of $a(p)/v(p)$.
    If this condition is not satisfied by the root $f_z(p^\ast,\mathbf{W}_V,\mathbf{W}_L) = 0$, the root is meaningless.
\end{rem}
\begin{thm}[Sufficient Condition for Solvability I]\label{suff_cond_solv_1}
    If the Riemann problem (\ref{mass_cons})-(\ref{mom_balance}), (\ref{init_data_rp}) is solvable without phase transition, see Subsection \ref{2p_wo_pt},
    then the same Riemann problem is also solvable taking into account phase transition according to the kinetic relation (\ref{kin_rel2}).
\end{thm}
\begin{prf}
    \newline
    \textbf{First Case $f(p^\ast,\mathbf{W}_V,\mathbf{W}_L) = 0$ with $p^\ast = p_0$:}
    In view of Section \ref{chap_pb_sol} we have $p_0 = \varphi(p^\ast)$, $z = 0$ and hence $f_z(p^\ast,\mathbf{W}_V,\mathbf{W}_L) = 0$.\\
    \newline
    \textbf{Second Case $f(p^\ast,\mathbf{W}_V,\mathbf{W}_L) = 0$ with $p^\ast > p_0$:} From that we have
    \begin{align*}
        \varphi(p^\ast) \stackrel{\ref{exis_uniq_thm}}{>} p^\ast > p_0\quad\text{and}\quad z(p^\ast) = p^\ast h(p^\ast,\varphi(p^\ast)) < 0.
    \end{align*}
    This gives
    \begin{align*}
        f_z(p^\ast,\mathbf{W}_V,\mathbf{W}_L) > f(p^\ast,\mathbf{W}_V,\mathbf{W}_L) = 0.
    \end{align*}
    So there exists a $p_V^\ast < p^\ast$ such that $f_z(p_V^\ast,\mathbf{W}_V,\mathbf{W}_L) = 0$.\\
    \newline
    \textbf{Third Case $f(p^\ast,\mathbf{W}_V,\mathbf{W}_L) = 0$ with $p^\ast < p_0$:} In this situation we obtain
    \begin{align*}
        0 = f(p^\ast,\mathbf{W}_V,\mathbf{W}_L) \stackrel{p^\ast < p_0}{<} f(p_0,\mathbf{W}_V,\mathbf{W}_L) \stackrel{\varphi(p_0) = p_0, z = 0}{=} f_z(p_0,\mathbf{W}_V,\mathbf{W}_L)
    \end{align*}
    Hence there exists a $p_V^\ast < p_0$ such that $f_z(p_V^\ast,\mathbf{W}_V,\mathbf{W}_L) = 0$.
\end{prf}
\begin{cor}
    Consider the Riemann problem (\ref{mass_cons})-(\ref{mom_balance}), (\ref{init_data_rp}) without phase transition and let $p^\ast$ be the solution for the pressure.
    Then we have for the same Riemann problem with phase transition and the corresponding solutions $p_V^\ast$ and $p_L^\ast = \varphi(p_V^\ast)$ the following relations:
    \begin{itemize}
        \item[(1)] $p^\ast = p_0$ implies $p_V^\ast = p_L^\ast = p_0$, i.e. equilibrium.
        \item[(2)] $p^\ast > p_0$ implies $p_0 < p_V^\ast < p^\ast$, i.e. condensation.
        \item[(3)] $p^\ast < p_0$ implies $p^\ast < p_L^\ast < p_0$, i.e. evaporation.
    \end{itemize}
\end{cor}
\begin{prf}
    The equilibrium case is obvious. The inequality $p_V^\ast < p^\ast$ in the second was obtained in the second part in the proof of Theorem \ref{suff_cond_solv_1}.
    It remains to show that $p_0 < p_V^\ast$. Assume that $p_V^\ast \leq p_0$, this gives
    \begin{align*}
        0 &= f_z(p_V^\ast,\mathbf{W}_V,\mathbf{W}_L) \leq f_z(p_0,\mathbf{W}_V,\mathbf{W}_L)\\
        &= f(p_0,\mathbf{W}_V,\mathbf{W}_L) \stackrel{p_0 < p^\ast}{<} f(p^\ast,\mathbf{W}_V,\mathbf{W}_L) = 0.
    \end{align*}
    For the evaporation case the inequality $p_L^\ast < p_0$ is a consequence of the third part in the proof of Theorem \ref{suff_cond_solv_1}.
    There we obtained $p_V^\ast < p_0$ and this gives, together with Theorem \ref{exis_uniq_thm}, the second inequality. Finally we want to prove the first inequality $p^\ast < p_L^\ast$.
    Again using Theorem \ref{exis_uniq_thm} gives $p_V^\ast > p_0$ if we assume $p_L^\ast > p_0$. By an analogous argument as for the second case this leads to a contradiction.
    Thus we have $p_L^\ast < p_0$.
\end{prf}
\begin{thm}[Sufficient Condition for Solvability II]\label{suff_cond_solv_2}
    Consider the Riemann problem (\ref{mass_cons})-(\ref{mom_balance}), (\ref{init_data_rp}) with phase transition.
    This Riemann problem is solvable by a (b)-type solution if and only if condition (\ref{vap_shock_z_bound}) holds and
    \begin{align*}
        f_z(\tilde{p},\mathbf{W}_V,\mathbf{W}_L) \geq 0.
    \end{align*}
\end{thm}
\begin{prf}
    The statement is obvious, since it guarantees a root for $f_z$.
\end{prf}
\newline
As in Subsection \ref{2p_wo_pt} we want to discuss wave patterns of type (a) and (c) for the Riemann problem (\ref{mass_cons})-(\ref{mom_balance}), (\ref{init_data_rp}) with phase transition.
The results are given in the subsequent three lemmata.
\begin{lem}\label{exclude_wavepattern_a}
    There is no solution with a wave pattern of type (a).
\end{lem}
\begin{prf}
    Assume there is a solution of type (a) as in Figure \ref{wave_pattern_a}.
    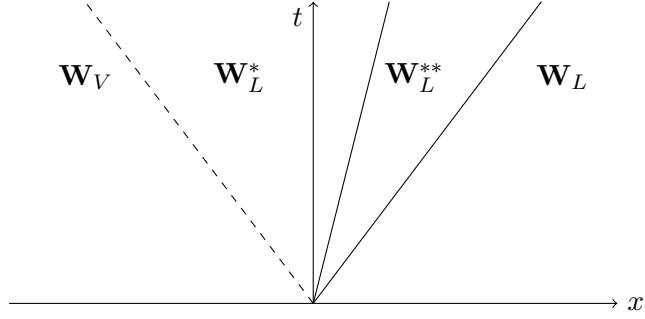
\begin{figure}[h!]
    \begin{center}
        \begin{tikzpicture}[domain=-4:4]
            %
            \draw[->] (-4,0) -- (4,0) node[right] {$x$};
            \draw[->] (0,0) -- (0,4);
            \draw[color=black] (0,3.8) node[left] {$t$};
            \draw (-2.5,3) node[left] {$\mathbf{W}_V$};
            \draw[color=black,dashed] (0,0) -- (-3,4);
            \draw (-0.5,3) node[left] {$\mathbf{W}_L^\ast$};
            \draw[color=black] (0,0) -- (1,4);
            \draw (1.85,3) node[left] {$\mathbf{W}_L^{\ast\ast}$};
            \draw (3.75,3) node[left] {$\mathbf{W}_L$};
            \draw[color=black] (0,0) -- (3,4);
        \end{tikzpicture}
    \end{center}
    \caption{Wave pattern of type (a) with notation}
    \label{wave_pattern_a}
    \end{figure}
    In this case we observe condensation and according to Corollary \ref{cond_evap_p_ineq} we have
    \begin{align*}
        z < 0\quad\text{and}\quad p_0 < p_V < p_L^\ast.
    \end{align*}
    Let us first assume that the left classical wave is a rarefaction wave. The head speed is given by $S = u_L^\ast - a_L(p_L^\ast)$ and we obtain
    \begin{align*}
        w = v_L(p_L^\ast)z + u_L^\ast \stackrel{(a)}{\leq} S = u_L^\ast - a_L(p_L^\ast)
        \quad\Leftrightarrow\quad
        z \leq -\frac{a_L(p_L^\ast)}{v_L(p_L^\ast)} \stackrel{(\ref{assumptions2})}{<} -\frac{a_V}{v_V}
    \end{align*}
    This is a contradiction and thus we can exclude this case. Given a shock instead of a rarefaction wave we have using (\ref{Q_square_ineq}) and the Lax condition
    \begin{align*}
        u_L^\ast - a_L(p_L^\ast) > S = u_L^\ast + v_L(p_L^\ast)Q_S > u_L^{\ast\ast} - a_L(p_L^{\ast\ast})
        \,\,\text{with}\,\,
        Q_S = -\frac{a_L(\bar{p}_L)}{v_L(\bar{p}_L)},\, \bar{p}_L \in (p_L^\ast,p_L^{\ast\ast}).
    \end{align*}
    Hence we yield
    \begin{align*}
        w < S
        \quad\Leftrightarrow\quad
        z < Q_S = -\frac{a_L(\bar{p}_L)}{v_L(\bar{p}_L)} \stackrel{p_L^\ast < \bar{p}_L}{<} -\frac{a_L(p_L^\ast)}{v_L(p_L^\ast)} \stackrel{(\ref{assumptions2})}{<} -\frac{a_V}{v_V}.
    \end{align*}
    Therefore we can also exclude this case and the proof is finished.
\end{prf}
\begin{lem}\label{2p_no_c_lem1}
    For the considered Riemann problem with phase transition exists no solution of type (c) with $p_L \geq p_0$.
\end{lem}
\begin{prf}
    A solution of type (c) implies an evaporation process which requires $p_L < p_0$.
\end{prf}
\begin{lem}\label{2p_no_c_lem2}
    For $p_L \in (\hat{p}_L,p_0]$ exists no solution of type (c) of the considered Riemann problem with phase transition.
\end{lem}
\begin{prf}
    Assume we have a wave pattern of type (c) as in Figure \ref{wave_pattern_c}. Hence we have evaporation and according to Corollary \ref{cond_evap_p_ineq} we have
    \begin{align*}
        z > 0\quad\text{and}\quad p_V^\ast < p_L < p_0.
    \end{align*}
    Let us first assume that the right classical wave is a rarefaction wave. The head speed is given by $S = u_V^\ast + a_V(p_V^\ast)$ and we obtain
    \begin{align*}
        w = v_V(p_V^\ast)z + u_V^\ast \stackrel{(c)}{\geq} S = u_V^\ast + a_V(p_V^\ast)
        \quad\Leftrightarrow\quad
        z \geq \frac{a_V(p_V^\ast)}{v_V(p_V^\ast)}.
    \end{align*}
    For a right shock ($Q_S > 0$) we have using (\ref{Q_square_ineq}) and the Lax condition
    \begin{align*}
        &u_V^{\ast\ast} + a_V(p_V^{\ast\ast}) > S = u_V^\ast + v_V(p_V^\ast)Q_S > u_V^\ast + a_V(p_V^\ast)\\
        &\text{with}\quad Q_S = \frac{a_V(\bar{p}_V)}{v_V(\bar{p}_V)},\quad \bar{p}_V \in (p_V^\ast,p_V^{\ast\ast}).
    \end{align*}
    Hence we yield
    \begin{align*}
        w > S
        \quad\Leftrightarrow\quad
        z > Q_S = \frac{a_V(\bar{p}_V)}{v_V(\bar{p}_V)} \stackrel{\bar{p}_V > p_V^\ast}{>} \frac{a_V(p_V^\ast)}{v_V(p_V^\ast)}.
    \end{align*}
    Due to Lemma \ref{z_max} we have an upper bound for the mass flux that does not initially exclude the conditions derived above for the rarefaction and shock wave.
    But the two cases are excluded if $z < a_V(p_V^\ast)/v_V(p_V^\ast)$. Indeed we have due to the monotonicity of $z$ and $a/v$ that
    \begin{align*}
        \exists\,\, \hat{p}_V < p_0\quad\text{such that}\quad \forall p_V \in (\hat{p}_V,p_0]: z(p_V) < \frac{a_V(p_V)}{v_V(p_V)}.
    \end{align*}
    Due to the strict monotonicity of $p_L^\ast = \varphi(p_V^\ast)$, see Theorem \ref{exis_uniq_thm}, the proof is complete.
\end{prf}

    \newpage
    \section{Phase Creation in Single Phase Flows}\label{chap_3p_rp}
\setcounter{para}{1}
\subsection{Condensation by Compression}\label{sec_nucl}
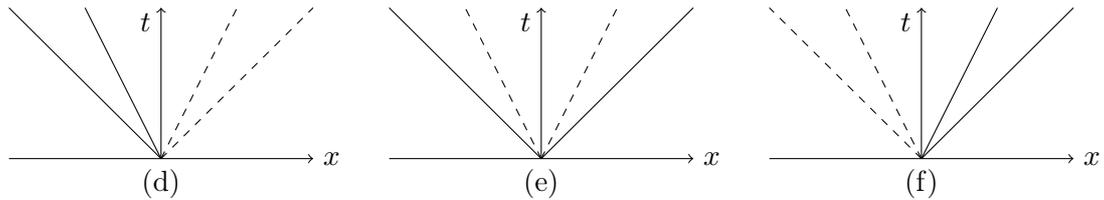
\begin{figure}[h!]
    \begin{center}
        \begin{tikzpicture}[domain=-7:7]
            %
            \draw[->] (-7,0) -- (-3,0) node[right] {$x$};
            \draw (-5,0) node[below] {(d)};
            \draw[->] (-5,0) -- (-5,2);
            \draw[color=black] (-5,1.8) node[left] {$t$};
            \draw[color=black] (-5,0) -- (-7,2);
            \draw[color=black] (-5,0) -- (-6,2);
            \draw[color=black,dashed] (-5,0) -- (-4,2);
            \draw[color=black,dashed] (-5,0) -- (-3,2);
            %
            \draw[->] (-2,0) -- (2,0) node[right] {$x$};
            \draw (0,0) node[below] {(e)};
            \draw[->] (0,0) -- (0,2);
            \draw[color=black] (0,1.8) node[left] {$t$};
            \draw[color=black] (0,0) -- (-2,2);
            \draw[color=black,dashed] (0,0) -- (-1,2);
            \draw[color=black,dashed] (0,0) -- (1,2);
            \draw[color=black] (0,0) -- (2,2);
            %
            \draw[->] (3,0) -- (7,0) node[right] {$x$};
            \draw (5,0) node[below] {(f)};
            \draw[->] (5,0) -- (5,2);
            \draw[color=black] (5,1.8) node[left] {$t$};
            \draw[color=black,dashed] (5,0) -- (3,2);
            \draw[color=black,dashed] (5,0) -- (4,2);
            \draw[color=black] (5,0) -- (6,2);
            \draw[color=black] (5,0) -- (7,2);
        \end{tikzpicture}
    \end{center}
    \caption{Wave patterns. Solid line: classical wave. Dashed line: phase boundary}
    \label{3_diff_wave_patterns2}
\end{figure}

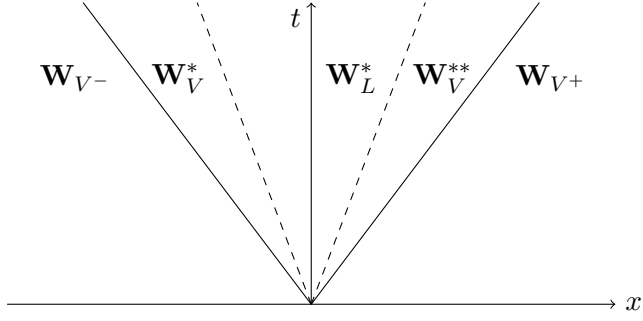
\begin{figure}[h!]
    \begin{center}
        \begin{tikzpicture}[domain=-4:4]
            %
            \draw[->] (-4,0) -- (4,0) node[right] {$x$};
            \draw[->] (0,0) -- (0,4);
            \draw[color=black] (0,3.8) node[left] {$t$};
            \draw (-2.5,3) node[left] {$\mathbf{W}_{V^-}$};
            \draw[color=black] (0,0) -- (-3,4);
            \draw (-1.25,3) node[left] {$\mathbf{W}_V^\ast$};
            \draw[color=black,dashed] (0,0) -- (-1.5,4);
            \draw (1,3) node[left] {$\mathbf{W}_L^\ast$};
            \draw[color=black,dashed] (0,0) -- (1.5,4);
            \draw (2.25,3) node[left] {$\mathbf{W}_V^{\ast\ast}$};
            \draw (3.75,3) node[left] {$\mathbf{W}_{V^+}$};
            \draw[color=black] (0,0) -- (3,4);
        \end{tikzpicture}
    \end{center}
    \caption{Wave pattern of type (e) with notation}
    \label{wave_pattern_e}
\end{figure}

Let us consider the following Riemann initial data with $\rho_{V^\pm} \in (0,\tilde{\rho}]$
\begin{align}
    \rho(x,0) = \begin{cases} &\rho_{V^-},\,x < 0\\ &\rho_{V^+},\,x > 0\end{cases}
    \quad\text{and}\quad
    u(x,0) = \begin{cases} &u_{V^-},\,x < 0\\ &u_{V^+},\,x > 0\end{cases}.\label{init_data_rp2}
\end{align}
Hence initially we have a Riemann problem for a single vapor phase and therefore we can directly apply the results obtained in subsection \ref{sec_riemann_problem}.
\begin{thm}[Solution of Isothermal Euler Equations for a Single Vapor Phase]\label{sol_rp_1p_vap}
    Let $f(p,\mathbf{W}_{V^-},\mathbf{W}_{V^+})$ be given as
    \begin{align*}
        f(p,\mathbf{W}_{V^-},\mathbf{W}_{V^+}) = f_-(p,\mathbf{W}_{V^-}) + f_+(p,\mathbf{W}_{V^+}) + \Delta u,\,\,\Delta u = u_{V^+} - u_{V^-},
    \end{align*}
    with the functions $f_-$ and$f_+$ given by
    \begin{align*}
        f_\pm(p,\mathbf{W}_{V^\pm}) &=
        \begin{dcases}
            \sqrt{-[\![p]\!][\![v]\!]},\,\,p > p_{V^\pm}\,\,\text{(Shock)}\\
            \int_{p_{V^\pm}}^p \frac{v_V(\zeta)}{a_V(\zeta)}\,\textup{d}\zeta,\,\,p \leq p_{V^\pm}\,\,\text{(Rarefaction)}
        \end{dcases}.
    \end{align*}
    If there is a root $f(p^\ast,\mathbf{W}_{V^-},\mathbf{W}_{V^+}) = 0$ with $0 < p^\ast \leq \tilde{p}$, then this root is unique.
    Further this is the unique solution for the pressure $p_V^\ast$ of the Riemann problem (\ref{mass_cons})-(\ref{mom_balance}), (\ref{init_data_rp2}).
    The velocity $u_V^\ast$ is given by
    \begin{align*}
        u_V^\ast = \frac{1}{2}(u_{V^+} + u_{V^-}) + \frac{1}{2}(f_+(p^\ast,\mathbf{W}_{V^+}) - f_-(p^\ast,\mathbf{W}_{V^-})).
    \end{align*}
\end{thm}
This is no new result and therefore it is well known, cf.\ Toro \cite{Toro2009}. Usually one looks for a pressure $p^\ast$ that solves $f(p,\mathbf{W}_{V^-},\mathbf{W}_{V^+}) = 0$.
Due to the asymptotic behavior there is always a solution. Nevertheless a solution with an unreasonable large vapor pressure is physically not meaningful,
since a sufficiently high pressure in a gas will lead to a phase transition to a liquid or even solid phase.
According to \cite{Hantke2013} we also only consider solutions which satisfy $0 < p^\ast \leq \tilde{p}$, where $\tilde{p}$ again denotes the maximal gas pressure.
This being said, we can find Riemann initial data without a solution. In this case proceed as follows.
\begin{defi}[Nucleation Criterion]\label{nucl_crit}
    If there is no solution of the Riemann problem (\ref{mass_cons})-(\ref{mom_balance}), (\ref{init_data_rp2}) according to Theorem \ref{sol_rp_1p_vap}, then nucleation occurs.
\end{defi}
If this criterion is fulfilled, we search a solution consisting of two classical waves and two phase boundaries. In the following we will again discuss several wave patterns.
\begin{lem}
    If there is a solution of the Riemann problem (\ref{mass_cons})-(\ref{mom_balance}), (\ref{init_data_rp2}) with two classical waves and two phase boundaries,
    then no wave is propagating inside the liquid phase. Hence classical waves may only occur in the vapor phase.
\end{lem}
\begin{prf}
    Assume a left going classical wave is propagating through the liquid phase. We denote the liquid states left and right to this wave by $\mathbf{W}_L^\ast$ and $\mathbf{W}_L^{\ast\ast}$.
    Further left to this classical wave there is a phase boundary moving with speed $w_1$. The vapor state left to this phase boundary is denoted by $\mathbf{W}_V^\ast$.
    Obviously this must be a condensation process and accordingly $p^\ast > p_0$ and $p_L^\ast > p_0$. This configuration is excluded by Lemma \ref{exclude_wavepattern_a}.
    Analogously we can discuss the case of a right going classical wave.
\end{prf}
\newline
As a consequence of the above result classical waves only propagate through the vapor phase. Hence we further have to discuss the following three patterns, see Figure \ref{3_diff_wave_patterns2}.
\begin{lem}\label{no_d_f_vap}
    There are no solutions of wave pattern types (d) and (f).
\end{lem}
\begin{prf}
    A solution with type (d) wave pattern corresponds to wave pattern type (c) in the previous Subsection \ref{2p_w_pt}, see Figure \ref{wave_pattern_c}.
    Thus by Lemma \ref{2p_no_c_lem1} and Lemma \ref{2p_no_c_lem2} we know that this is only possible for sufficiently small pressures and therefore implies evaporation.
    Since we have a condensation process wave pattern type (d) can be excluded. Analogously we discuss a type (f) solution.
    This corresponds to a type (a) solution in Subsection \ref{2p_w_pt}, see Figure \ref{wave_pattern_a}.
    Hence due to Lemma \ref{exclude_wavepattern_a} a solution of wave pattern type (f) is also impossible.
\end{prf}
\newline
Consequently the only possible wave pattern in this case is of type (e), see Figure \ref{wave_pattern_e}.
\begin{lem}
    For a solution of wave pattern type (e) the equality $p_V^\ast = p_V^{\ast\ast}$ holds.
\end{lem}
\begin{prf}
    Across the left phase boundary the liquid pressure $p_L^\ast$ is uniquely defined by the vapor pressure $p_V^\ast$ using Theorem \ref{exis_uniq_thm}.
    So far we assumed the vapor left of the liquid phase. For the right phase boundary the opposite is the case and thus we have to use the kinetic relation (\ref{kin_rel_lv}).
    Nevertheless the results of the previous section remain unchanged and hence we obtain the same function to determine the liquid pressure
    \begin{align*}
        p_L^\ast = \varphi(p_V^\ast) = \varphi(p_V^{\ast\ast}).
    \end{align*}
    Hence the vapor pressures are equal.
\end{prf}
\newline
Taking into account that there are two phase boundaries and using the results obtained in the previous sections we can state the following theorem.
\begin{thm}[Solution of Isothermal Euler Equations for Two Vapor States with Phase Transition]\label{cond_compr_thm}
    Consider the Riemann problem (\ref{mass_cons})-(\ref{mom_balance}), (\ref{init_data_rp2}) and assume the nucleation criterion is satisfied.
    Let $f_z(p,\mathbf{W}_{V^-},\mathbf{W}_{V^+})$ be given as
    \begin{align*}
        f_z(p,\mathbf{W}_{V^-},\mathbf{W}_{V^+}) = f_-(p,\mathbf{W}_{V^-}) + f_+(p,\mathbf{W}_{V^+}) + 2z[\![v]\!] + \Delta u,\,\,\Delta u = u_{V^+} - u_{V^-},
    \end{align*}
    with the functions $f_-$ and$f_+$ given by
    \begin{align*}
        f_\pm(p,\mathbf{W}_{V^\pm}) &=
        \begin{dcases}
            \sqrt{-[\![p]\!][\![v]\!]},\,\,p > p_{V^\pm}\,\,\text{(Shock)}\\
            \int_{p_{V^\pm}}^p \frac{v_V(\zeta)}{a_V(\zeta)}\,\textup{d}\zeta,\,\,p \leq p_{V^\pm}\,\,\text{(Rarefaction)}
        \end{dcases}.
    \end{align*}
    Here $z$ is given by (\ref{kin_rel1}) and $[\![v]\!] = v_L(\varphi(p)) - v_V(p)$. The function $\varphi$ is implicitly defined by (\ref{pb_mom_balance}).
    If there is a root $f_z(p^\ast,\mathbf{W}_{V^-},\mathbf{W}_{V^+}) = 0$ with $p_0 < p^\ast \leq \tilde{p}$, then this root is the only one.
    Further this is the unique solution for the vapor pressures $p_V^\ast = p_V^{\ast\ast}$ of the Riemann problem (\ref{mass_cons})-(\ref{mom_balance}), (\ref{init_data_rp2}).
    The liquid velocity $u_L^\ast$ is given by
    \begin{align*}
        u_L^\ast = \frac{1}{2}(u_{V^+} + u_{V^-}) + \frac{1}{2}(f_+(p^\ast,\mathbf{W}_{V^+}) - f_-(p^\ast,\mathbf{W}_{V^-})).
    \end{align*}
\end{thm}
By the previous results it is obvious that $f_z$ has at most one root. By construction this root is the solution for the vapor pressure in the two star regions, see Figure \ref{wave_pattern_e}.
The following results are completely analogue to those obtained before for the two phase case.
\begin{rem}
    Note that $u_V^\ast \neq u_V^{\ast\ast}$ with $u_V^\ast + u_V^{\ast\ast} = 2u_L^\ast$.
\end{rem}
\begin{thm}[Sufficient Condition for Solvability I]\label{suff_cond_sol3}
    Consider the Riemann problem (\ref{mass_cons})-(\ref{mom_balance}), (\ref{init_data_rp2}). This problem is solvable without phase transition if and only if
    \begin{align*}
        f(\tilde{p},\mathbf{W}_{V^-},\mathbf{W}_{V^+}) \geq 0.
    \end{align*}
    Here $f$ is given as in Theorem \ref{sol_rp_1p_vap}.
\end{thm}
\begin{prf}
    This statement is obvious due to the monotonicity of $f$.
\end{prf}
\begin{thm}[Sufficient Condition for Solvability II]\label{suff_cond_sol4}
    Consider the Riemann problem (\ref{mass_cons})-(\ref{mom_balance}), (\ref{init_data_rp2}) and assume the nucleation criterion is satisfied.
    Accounting for phase transition, this problem is solvable if and only if
    \begin{align*}
        f_z(\tilde{p},\mathbf{W}_{V^-},\mathbf{W}_{V^+}) \geq 0.
    \end{align*}
    The function $f_z$ is defined as in Theorem \ref{cond_compr_thm}.
\end{thm}
\begin{prf}
    Again the statement is obvious due to the monotonicity of $f_z$.
\end{prf}
\subsection{Evaporation by Expansion}\label{sec_cav}
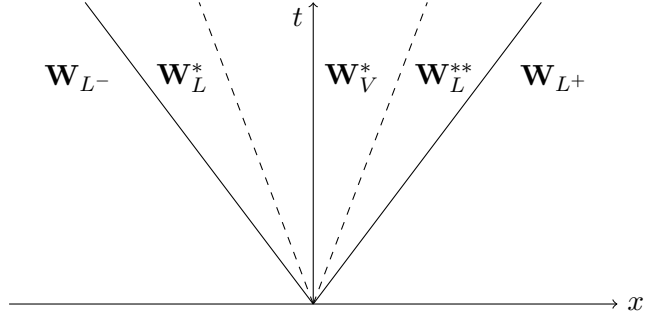
\begin{figure}[h!]
    \begin{center}
        \begin{tikzpicture}[domain=-4:4]
            %
            \draw[->] (-4,0) -- (4,0) node[right] {$x$};
            \draw[->] (0,0) -- (0,4);
            \draw[color=black] (0,3.8) node[left] {$t$};
            \draw (-2.5,3) node[left] {$\mathbf{W}_{L^-}$};
            \draw[color=black] (0,0) -- (-3,4);
            \draw (-1.25,3) node[left] {$\mathbf{W}_L^\ast$};
            \draw[color=black,dashed] (0,0) -- (-1.5,4);
            \draw (1,3) node[left] {$\mathbf{W}_V^\ast$};
            \draw[color=black,dashed] (0,0) -- (1.5,4);
            \draw (2.25,3) node[left] {$\mathbf{W}_L^{\ast\ast}$};
            \draw (3.75,3) node[left] {$\mathbf{W}_{L^+}$};
            \draw[color=black] (0,0) -- (3,4);
        \end{tikzpicture}
    \end{center}
    \caption{Wave pattern of type (e) with notation for the liquid case}
    \label{wave_pattern_e2}
\end{figure}

Now we consider the following Riemann initial data with $\rho_{L^\pm} \geq \rho_L^{min}$
\begin{align}
    \rho(x,0) = \begin{cases} &\rho_{L^-},\,x < 0\\ &\rho_{L^+},\,x > 0\end{cases}
    \quad\text{and}\quad
    u(x,0) = \begin{cases} &u_{L^-},\,x < 0\\ &u_{L^+},\,x > 0\end{cases}.\label{init_data_rp3}
\end{align}
Hence initially we have a Riemann problem for a single liquid phase. We have seen so far that at a planar phase boundary the liquid pressure is always positive.
However it is known that negative liquid pressures are possible, cf.\,\,Davitt et al.\,\,\cite{Davitt2010} for water. This gives rise to cavitation in the liquid phase. 
Again, in the liquid-vapor case a negative liquid pressure is forbidden, see (\ref{p_jump_pos}). Nevertheless in the liquid-liquid case we may encounter negative liquid pressures.
We define the smallest possible liquid pressure to be $p_{min}$ and with this definition we obtain the following result.
\begin{thm}[Solution of Isothermal Euler Equations for a Single Liquid Phase]\label{sol_rp_1p_liq}
    Let $f(p,\mathbf{W}_{L^-},\mathbf{W}_{L^+})$ be given as
    \begin{align*}
        f(p,\mathbf{W}_{L^-},\mathbf{W}_{L^+}) = f_-(p,\mathbf{W}_{L^-}) + f_+(p,\mathbf{W}_{L^+}) + \Delta u,\,\,\Delta u = u_{L^+} - u_{L^-},
    \end{align*}
    with the functions $f_-$ and $f_+$ given by
    \begin{align*}
        f_\pm(p,\mathbf{W}_{L^\pm}) &=
        \begin{dcases}
            \sqrt{-[\![p]\!][\![v]\!]},\,\,p > p_{L^\pm}\,\,\text{(Shock)}\\
            \int_{p_{L^\pm}}^p \frac{v_L(\zeta)}{a_L(\zeta)}\,\textup{d}\zeta,\,\,p \leq p_{L^\pm}\,\,\text{(Rarefaction)}
        \end{dcases}.
    \end{align*}
    If there is a root $f(p^\ast,\mathbf{W}_{L^-},\mathbf{W}_{L^+}) = 0$ with $p_{min} \leq p^\ast$, then this root is unique.
    Further this is the unique solution for the pressure $p_L^\ast$ of the Riemann problem (\ref{mass_cons})-(\ref{mom_balance}), (\ref{init_data_rp3}).
    The velocity $u_L^\ast$ is given by
    \begin{align*}
        u_L^\ast = \frac{1}{2}(u_{L^+} + u_{L^-}) + \frac{1}{2}(f_+(p^\ast,\mathbf{W}_{L^+}) - f_-(p^\ast,\mathbf{W}_{L^-})).
    \end{align*}
\end{thm}
\begin{rem}
    For simplicity we choose $p_{min} = 0$ but lower values are possible and the theoretical results do not depend on the specific value of $p_{min}$.
\end{rem}
Analogously to the case of nucleation we define the following.
\begin{defi}[Cavitation Criterion]\label{cav_crit}
    If there is no solution of the Riemann problem (\ref{mass_cons})-(\ref{mom_balance}), (\ref{init_data_rp3}) according to Theorem \ref{sol_rp_1p_liq}, then cavitation may occur.
\end{defi}
If this criterion is fulfilled, we look for a solution involving a vapor phase, two phase boundaries and two classical waves. Again we discuss the possible patterns.
\begin{lem}
    Assume there is a solution of the Riemann problem (\ref{mass_cons})-(\ref{mom_balance}), (\ref{init_data_rp3}) consisting of two classical waves and two phase boundaries.
    If the pressures $p_{L^-}, p_{L^+}$ are sufficiently large then no wave travels through the vapor.
\end{lem}
The proof is analogue to the one of Lemma \ref{2p_no_c_lem2}.
\begin{lem}
    There is no solution of type (d) or (f); see Figure \ref{3_diff_wave_patterns2}.
\end{lem}
The proof of this lemma is analogue to the one of Lemma \ref{no_d_f_vap}. Accordingly we construct solutions of type (e), see Figure \ref{wave_pattern_e2}.
\begin{thm}[Solution of Isothermal Euler Equations for Two Liquid States with Phase Transition]\label{evap_exp_thm}
    Consider the Riemann problem (\ref{mass_cons})-(\ref{mom_balance}), (\ref{init_data_rp3}) and assume the cavitation criterion is satisfied.
    Let $f_z(p,\mathbf{W}_{L^-},\mathbf{W}_{L^+})$ be given as
    \begin{align*}
        f_z(p,\mathbf{W}_{L^-},\mathbf{W}_{L^+}) = f_-(p,\mathbf{W}_{L^-}) + f_+(p,\mathbf{W}_{L^+}) + 2z[\![v]\!] + \Delta u,\,\,\Delta u = u_{L^+} - u_{L^-},
    \end{align*}
    with the functions $f_-$ and $f_+$ given by
    \begin{align*}
        f_\pm(p,\mathbf{W}_{L^\pm}) &=
        \begin{dcases}
            \sqrt{-[\![p]\!][\![v]\!]},\,\,\varphi(p) > p_{L^\pm}\,\,\text{(Shock)}\\
            \int_{p_{L^\pm}}^{\varphi(p)} \frac{v_L(\zeta)}{a_L(\zeta)}\,\textup{d}\zeta,\,\,\varphi(p) \leq p_{L^\pm}\,\,\text{(Rarefaction)}
        \end{dcases}.
    \end{align*}
    Here $z$ is given by (\ref{kin_rel1}) and $[\![v]\!] = v_L(\varphi(p)) - v_V(p)$. The function $\varphi$ is implicitly defined by (\ref{pb_mom_balance}).
    If there is a root $f_z(p^\ast,\mathbf{W}_{L^-},\mathbf{W}_{L^+}) = 0$ with $p_{min} \leq p^\ast$, then this root is unique.
    Further this is the unique solution for the vapor pressures $p_V^\ast$ of the Riemann problem (\ref{mass_cons})-(\ref{mom_balance}), (\ref{init_data_rp3}).
    The vapor velocity $u_V^\ast$ is given by
    \begin{align*}
        u_V^\ast = \frac{1}{2}(u_{L^+} + u_{L^-}) + \frac{1}{2}(f_+(p^\ast,\mathbf{W}_{L^+}) - f_-(p^\ast,\mathbf{W}_{L^-})).
    \end{align*}
\end{thm}
\begin{prf}
    Due to the previous results, the function $f_z$ has at most one root. This root is by construction the solution for the vapor pressure in the star region.
\end{prf}
\newline
Completely analogue to the condensation case, see Subsection \ref{sec_nucl}, we have the following results.
\begin{thm}[Sufficient Condition for Solvability I]\label{suff_cond_sol5}
    Consider the Riemann problem (\ref{mass_cons})-(\ref{mom_balance}), (\ref{init_data_rp3}). This problem is solvable without phase transition if and only if
    \begin{align*}
        f(p_{min},\mathbf{W}_{L^-},\mathbf{W}_{L^+}) \leq 0.
    \end{align*}
    Here $f$ is given as in Theorem \ref{sol_rp_1p_liq}.
\end{thm}
\begin{prf}
    The statement is easily verified due to the monotonicity of $f$.
\end{prf}
\begin{thm}[Sufficient Condition for Solvability II]\label{suff_cond_sol6}
    Consider the Riemann problem (\ref{mass_cons})-(\ref{mom_balance}), (\ref{init_data_rp3}) and assume the cavitation criterion is satisfied.
    If we admit phase transition, this problem is always solvable.
\end{thm}
\begin{prf}
    This statement is obvious due to the fact that $z[\![v]\!] \to -\infty$ for $p_V^\ast \to 0$.
\end{prf}

    \newpage
    \section{Conclusion}\label{chap_conclusion}
\setcounter{para}{1}
\subsection{Discussion of the Assumptions}\label{discuss_assumptions}
In this part we now want to discuss the assumptions previously made to solve the problem.
Basically we have three types of requirements. First there are the ones due to the underlying thermodynamics, in particular the first and second law of thermodynamics.
Second there are conditions, one needs to solve the single phase Riemann problem for the Euler equations.
The third type concerns the assumptions imposed to solve the two phase problem.
Note that the assumptions are sufficient, from a mathematical point of view, to obtain the results presented throughout this work.\\
\newline
From a thermodynamic point of view we have first and foremost to satisfy the first and second law of thermodynamics including the requirement of thermodynamic stability (\ref{thermodyn_stab}).
This is obtained by deriving the pressure law from a suited thermodynamic potential.\\
The conditions imposed on the EOS in order to solve the (single phase) Riemann problem for the Euler equations are
\begin{align*}
    \gamma > 0,\quad\mathcal{G} > 0,\quad v(p) \stackrel{p\to\infty}\to 0,\quad\text{and}\quad v(p) \stackrel{p \to 0}{\to} \infty
\end{align*}
That we require the single phase Riemann problem to be solvable is of course reasonable, since otherwise any further discussion would be unnecessary.
The conditions above are completely analogue to those stated in \cite{Menikoff1989}.
Note that for any EOS where the speed of sound is a constant (such as in \cite{Hantke2013}) we have $\mathcal{G} = 1$.
We want to point out that the aforementioned requirements of type one and two are basically no new or additional assumptions since they are already needed to treat the single phase case.\\
\newline
Since we are concerned with discussing the case of two phases it is reasonable to assume that all single phase requirements are met and only a few new ones need to be added.
In order to solve the two phase Riemann problem we need the additional assumptions (\ref{assumptions1}) and (\ref{assumptions2}).\\
The uniform upper bound for the quotient of the specific volumes basically tells us how close we can get to the critical point, where the volumes would become equal.
The case of $\alpha = 1$, i.e. we include the critical point where the volumes become equal, is not considered here and has to be treated separately.\\
The constant $\beta$ bounds the quotient of the sound speeds and is only needed to be strict smaller than $1/\alpha$.\\
The assumption on the lower bound of $\gamma_V$ in (\ref{assumptions1}) is a rather technical one.
Nevertheless if we assume $\tau$ to be as in (\ref{tau_special}) and consider the ideal gas EOS for the vapor phase we have
\begin{align*}
    1 = \gamma_V > \tau a_V^3(1 - \alpha)^2 =  \frac{1}{\sqrt{2\pi}}\left(\frac{m}{kT_0}\right)^{\frac{3}{2}}(1 - \alpha)^2a_V^3 = \frac{(1 - \alpha)^2}{\sqrt{2\pi}}.
\end{align*}
Hence this bound is easily satisfied. If the sound speed of the vapor phase would depend on the pressure one would have to check this requirement more carefully.
We further want to emphasize that apart from $\tau > 0$ and $(\ref{assumptions1})_3$ we do not assume any particular shape or even magnitude of $\tau$.
This is a further key point that contributes to the generality of our result.\\
The last requirement in (\ref{assumptions1}) is concerned with the maximum vapor pressure.
Due to this bound the vapor is allowed to be compressed (depending on $\alpha$) with more than the saturation pressure. This allows metastable states, which is reflected in the Maxwell construction.
Here of course one has to make sure that the maximum vapor pressure $\tilde{p}$ defined in Definition \ref{max_pV} satisfies this bound.
This can be guaranteed by choosing an appropriate temperature and also how the two EOS are connected by $\bar{v}(p)$ in Definition \ref{max_pV}.
Usually $\tilde{p}$ will only be slightly larger than the saturation pressure for a wide range of temperatures.\\
\newline
Now we want to comment assumptions (\ref{assumptions2}). Let us first consider $\gamma_L$.
Over wide temperature ranges we have $\gamma_L \geq 1$ for many substances.
For example in Section \ref{examples} we consider the \textit{linear} and \textit{nonlinear Tait} EOS for liquid water and for this type of EOS modeling water this is true up to $636.165 K$.
A similar result can be obtained for the \textit{van der Waals} EOS.
Above that temperature it is not possible to use the ideal (polytropic) gas EOS together with such a liquid EOS, because it would contradict (\ref{assumptions2}) (\textbf{II}).\\
Regarding case (\textbf{II}) we want to emphasize that for $1 > \gamma_L > \gamma_V$ the inequality including $\alpha$ is trivial.
In fact in numerical studies we exemplary obtained that this property is also true for the van der Waals EOS up to $\approx 640 K$.\\
Now we want to comment on $\varepsilon(\gamma_V)$ in (\ref{assum_const}). Using the ideal gas or the polytropic gas EOS gives $\gamma_V = 1$ and hence
\begin{align*}
    \varepsilon_0 := \varepsilon(1) = - \tau a_V^3(1 - \alpha)^2\left(1 - (\alpha\beta)^2\right) < 0.
\end{align*}
We consider (\ref{assumptions2}) (\textbf{II}) and have $1 + \varepsilon_0/\alpha < 0$ over large temperature ranges.
Suppose this term becomes positive at high temperatures, it is however still smaller than one. Whereas at the same time $\gamma_L$ approaches one. Hence this bound may be still valid.
This of course has to be checked for any given EOS.
\subsection{Examples}\label{examples}
Now we want to present several examples of choices for the equations of state used to model the fluid under consideration, in this case water.
First we will discuss the ideal gas EOS for the vapor phase and the (linear) Tait EOS for the liquid phase as in \cite{Hantke2013}. For the ideal gas we obtain
\begin{align}
    p_V(v_V) = \frac{kT_0}{m}\frac{1}{v_V},\quad \gamma_V = 1,\quad \mathcal{G}_V = 1.\label{id_gas_prop1}
\end{align}
Here $k$ is the Boltzmann constant, $T_0$ is the fixed temperature and $m$ denotes the mass of a single water molecule. Considering the liquid phase we obtain
\begin{align}
    p_L(v_L) = p_0 + K_0\left(\frac{v_0}{v_L} - 1\right),\quad\gamma_L = \left(1 + \frac{v_L}{K_0v_0}(p_0 - K_0)\right)^{-1}\stackrel{K_0 \geq p_0}{\geq} 1,\quad\mathcal{G}_L = 1.\label{lin_tait_prop}
\end{align}
The quantities with index zero are calculated at the saturation state corresponding to $T_0$.
We further have the saturation pressure $p_0$, the modulus of compression $K_0$ and the specific liquid volume $v_0$, cf. \cite{Wagner1998}.
Note that the relation $K_0 \geq p_0$ breaks down for temperatures \textit{above} $636.165 K$ ($T_c = 647.096 K$).\\
Both EOS are linear functions of the mass density and thus it is reasonable to connect them with a linear function $\bar{p}$.
Hence we obtain the specific volume of the vapor phase corresponding to the maximum vapor pressure $\tilde{p}$ according to Definition \ref{max_pV} as the solution of the following equation
\begin{align}
    0 = K_0v_0\ln\frac{v_0}{v_m} + \frac{\tilde{v}}{v_m - \tilde{v}}\frac{kT_0}{m}\ln\frac{\tilde{v}}{v_m} + \frac{kT_0}{m}\ln\frac{v_V(p_0)}{v_m}.\label{maxw_constr1}
\end{align}
Here $v_{m}$ is chosen such that
\begin{align*}
    %
    v_{m} = \begin{cases}
                &v_0\left(1 - \dfrac{p_0}{K_0}\right)^{-1}, T_0 \leq 620 K,\\
                &v_0\left(1 + \dfrac{T_c - T_0}{T_c}\right), T_0 > 620 K.
            \end{cases}
\end{align*}
Using (\ref{maxw_constr1}) we can calculate the quotient $v_{m}/\tilde{v}$ for every reasonable temperature and thus obtain $\alpha$ and also $\beta$.
Now we can check the assumptions given in (\ref{assumptions1}), (\ref{assumptions2}). We have for temperatures up to $636.165 K$ the following
\begin{align*}
    \alpha \lesssim 0.1949,\quad \alpha\beta \lesssim 0.5419,\quad \frac{1}{\gamma_L} - \left(1 + \frac{\varepsilon_0}{\alpha}\right) \gtrsim 0.7484 \quad\text{and}\quad \tilde{p} \lesssim 1.4825p_0.
\end{align*}
Thus all requirements are met as expected and the limiting factor here are not the assumptions but the choice of the EOS.
\begin{rem}
    Note that in the isothermal case the linear Tait EOS is equivalent to the \textit{stiffend gas} EOS
    \begin{align*}
        p_L(v_L) = C(\gamma - 1)\frac{T_0}{v_L} - p_c\quad\text{with}\quad C(\gamma - 1) = \frac{K_0v_0}{T_0}\quad\text{and}\quad p_c = K_0 - p_0.
    \end{align*}
\end{rem}
As a second example we want to use the nonlinear Tait EOS instead of the linear one, i.e.
\begin{align}
    p_L(v_L) = p_0 + K_0\left(\left(\frac{v_0}{v_L}\right)^\nu - 1\right),\,\,\nu > 1.\label{nonlin_tait_eos}
\end{align}
Again we use the ideal gas EOS for the vapor phase. We obtain for the nonlinear Tait EOS
\begin{align*}
    \gamma_L = \nu\left(1 + \left(\frac{v_L}{v_0}\right)^\nu\left(\frac{p_0}{K_0} - 1\right)\right)^{-1} > 1,\quad \mathcal{G}_L = \frac{\nu + 1}{2}
    %
\end{align*}
and
\begin{align*}
    %
    v_{m} = \begin{cases}
                &v_0\left(1 - \dfrac{p_0}{K_0}\right)^{-\frac{1}{\nu}}, T_0 \leq 620 K,\\
                &v_0\left(1 + \dfrac{T_c - T_0}{T_c}\right), T_0 > 620 K.
            \end{cases}
\end{align*}
Next with an approach analogue to the previous case we obtain $\tilde{v}$ as solution of the following equation and then calculate $\tilde{p}$
\begin{align}
    0 &= (p_0 - K_0)v_{m} + K_0v_0 + \frac{K_0v_0^\nu}{1 - \nu}\left(\frac{1}{v_{m}^{\nu - 1}} - \frac{1}{v_0^{\nu - 1}}\right)\notag\\
      &+ \frac{v_{m}}{v_{m} - \tilde{v}}\frac{kT_0}{m}\ln\frac{\tilde{v}}{v_{m}} + \frac{kT_0}{m}\ln\frac{v_V(p_0)}{\tilde{v}}.\label{maxw_constr2}
\end{align}
We can use (\ref{maxw_constr2}) to calculate the quotient $v_{m}/\tilde{v}$ for every reasonable temperature and thus obtain $\alpha$ and also $\beta$.
Here we use $\nu = 7$ as in \cite{Saurel1999}.
We again check the assumptions given in (\ref{assumptions1}), (\ref{assumptions2}) and obtain for temperatures up to $636.165 K$ the following
\begin{align*}
    \alpha \lesssim 0.1645,\quad \alpha\beta \lesssim 0.1818,\quad \frac{1}{\gamma_L} - \left(1 + \frac{\varepsilon_0}{\alpha}\right) \gtrsim 0.7795 \quad\text{and}\quad \tilde{p} \lesssim 1.2511p_0.
\end{align*}
Hence this choice of EOS is also suitable for solving this problem for temperatures from $273.15 K$ up to $636.165 K$.
Again the limiting factor here are not the assumptions but the choice of the EOS.
\subsection{Conclusion}\label{conclusion}
The aim of the present work was to investigate the Riemann problem for the isothermal Euler equations when liquid and vapor phases are present which may condensate or evaporate.
We proved that there exist unique solutions under the given assumptions.
To this end we allow any EOS which satisfies these assumptions, especially nonlinear ones.
This is a huge improvement to the previous work \cite{Hantke2013} where only two specific linear EOS were chosen to solve this problem.
In contrast to this we for example allow the speed of sound to depend on the pressure or volume instead of being constant.
Furthermore allow phase transitions where the pressures are not in equilibrium as for example in \cite{Dumbser2013}. Additionally we can treat nucleation an cavitation.
In view of the work by Hantke, T. \cite{Hantke2015a} the last point has to be emphasized.
To our knowledge this is the most general result concerning Riemann problems for isothermal two phase flows.

    \cleardoublepage

    \phantomsection
    \bibliographystyle{abbrv}
    \bibliography{literatur_2p_iso_gen}

\end{document}